\documentclass[11pt,letterpaper,reqno]{amsart}

\title[Nonlinear Stability of Planar $N$-Vortex Problem]{Nonlinear Stability of Relative Equilibria in the Planar $N$-Vortex Problem with Non-Zero Total Circulation}
\author{Tomoki Ohsawa}
\address{Department of Mathematical Sciences, The University of Texas at Dallas, 800 W Campbell Rd, Richardson, TX 75080-3021}
\email{tomoki@utdallas.edu}
\date{\today}

\keywords{$N$-vortex problem, stability of relative equilibria, Lie--Poisson equation}

\subjclass[2020]{37J25,53D20,70G65,70H14,76B47}

\usepackage{graphicx,amssymb,paralist}

\usepackage{caption,subcaption}
\usepackage[margin=1in, marginpar=.5in]{geometry}

\usepackage[numbers,sort&compress]{natbib}

\usepackage[colorlinks=false]{hyperref}

\usepackage[nameinlink]{cleveref}

\usepackage{tikz}
\usetikzlibrary{tikzmark,fit}
\usepackage{eso-pic}
\usepackage{tikz-cd}

\usepackage{nicematrix}

\theoremstyle{plain}
\newtheorem{theorem}{Theorem}[section]

\newtheorem{lemma}[theorem]{Lemma}
\newtheorem{proposition}[theorem]{Proposition}

\theoremstyle{definition}

\theoremstyle{remark}
\newtheorem{remark}[theorem]{Remark}

\def\od#1#2{\dfrac{d#1}{d#2}}

\def\fd#1#2{\frac{\delta #1}{\delta #2}}

\def\tpd#1#2{\partial #1/\partial #2}

\def\parentheses#1{\!\left(#1\right)}
\def\brackets#1{\!\left[#1\right]}
\def\braces#1{\!\left\{#1\right\}}

\def\tr{\mathop{\mathrm{tr}}\nolimits}
\def\rank{\operatorname{rank}}

\def\norm#1{\left\|#1\right\|}

\def\DS{\displaystyle}
\def\R{\mathbb{R}}
\def\C{\mathbb{C}}

\def\N{\mathbb{N}}
\def\defeq{\mathrel{\mathop:}=}
\def\eqdef{=\mathrel{\mathop:}}
\def\setdef#1#2{ \left\{ #1 \,\mid\, #2 \right\} }
\def\ip#1#2{\left\langle#1,#2\right\rangle}

\renewcommand{\Re}{\operatorname{Re}}
\renewcommand{\Im}{\operatorname{Im}}
\def\eps{\varepsilon}
\def\rmi{{\rm i}}

\def\d{\mathbf{d}}
\def\ins#1{{\bf i}_{#1}}
\def\PB#1#2{\left\{#1,#2\right\}}

\newcommand\ad{\operatorname{ad}}

\def\SO{\mathsf{SO}}
\def\SE{\mathsf{SE}}

\def\u{\mathfrak{u}}

\def\vK{\mathfrak{v}_{\mathcal{K}}}
\def\nonzerovK{\mathring{\mathfrak{v}}_{\mathcal{K}}}

\begin{document}

\footskip=.6in

\begin{abstract}
  We prove a sufficient condition for nonlinear stability of relative equilibria in the planar $N$-vortex problem with non-zero total circulation; the stability is in the sense of Lyapunov to perturbations that preserve the momentum map.
  Using the condition, we also prove nonlinear stability of some known relative equilibria beyond the ranges of parameters that were previously known to be stable.
  Our result builds on our previous work on the Hamiltonian formulation of its relative dynamics as a Lie--Poisson system.
  The relative dynamics recasts the relative equilibria of the $N$-vortex problem as fixed points in the Lie--Poisson relative dynamics.
  We analyze the stability of such fixed points by exploiting the Hamiltonian formulation as well as invariants and constraints that naturally arise in the relative dynamics.
  We apply the method to the rhombus equilibria with two vortices of circulation 1 and the other two of circulation $\gamma$, and prove that they are stable for $-2 + \sqrt{3} < \gamma < 0$, which were previously known to be only linearly stable due to Roberts.
\end{abstract}
  
\maketitle

\section{Introduction}
\subsection{Relative Dynamics of $N$ Point Vortices}
Let $\{ \mathbf{x}_{i} = (x_{i},y_{i}) \in \R^{2} \}_{i=1}^{N}$ be the positions on the plane $\R^{2}$ of $N$ point vortices (vortices for short in what follows) with circulations $\{ \Gamma_{i} \in \R\backslash\{0\} \}_{i=1}^{N}$.
The dynamics of these vortices is governed by the following system of equations (see, e.g., \citet[Section~2.1]{Ne2001} and \citet[Section~2.1]{ChMa1993}):
For every $i \in \{1, \dots, N\}$,
\begin{equation}
  \label{eq:vortices}
  \dot{x}_{i} = -\frac{1}{2\pi} \sum_{\substack{1\le j \le N\\ j \neq i}} \Gamma_{j} \frac{y_{i} - y_{j}}{\norm{\mathbf{x}_{i} - \mathbf{x}_{j}}^{2}},
  \qquad
  \dot{y}_{i} = \frac{1}{2\pi} \sum_{\substack{1\le j \le N\\ j \neq i}} \Gamma_{j} \frac{x_{i} - x_{j}}{\norm{\mathbf{x}_{i} - \mathbf{x}_{j}}^{2}},
  \qquad
\end{equation}
or via the identification of $\R^{2}$ with $\C$ by $(x_{i},y_{i}) \mapsto x_{i} + \rmi y_{i} \eqdef q_{i}$,
\begin{equation}
  \label{eq:vortices-C}
  \dot{q}_{i} = \frac{\rmi}{2\pi} \sum_{\substack{1\le j \le N\\ j \neq i}} \Gamma_{j} \frac{q_{i} - q_{j}}{|q_{i} - q_{j}|^{2}}.
\end{equation}

It is well known that the above system is a Hamiltonian system on $\R^{2N} = \{ (\mathbf{x}_{1}, \dots, \mathbf{x}_{N}) \}$ or $\C^{N} = \{ (q_{1}, \dots, q_{N}) \}$ with the symplectic form
\begin{equation}
  \label{eq:Omega}
  \Omega
  \defeq \sum_{i=1}^{N}\Gamma_{i} \d{x}_{i} \wedge \d{y}_{i}
  = -\frac{1}{2} \sum_{i=1}^{N} \Gamma_{i} \Im(\d{q}_{i} \wedge \d{q}_{i}^{*})
\end{equation}
and the Hamiltonian
\begin{align*}
  H(\mathbf{x}_{1}, \dots, \mathbf{x}_{N})
  &\defeq -\frac{1}{4\pi} \sum_{1\le i < j \le N} \Gamma_{i} \Gamma_{j} \ln\norm{\mathbf{x}_{i} - \mathbf{x}_{j}}^{2} \\
  &= -\frac{1}{4\pi} \sum_{1\le i < j \le N} \Gamma_{i} \Gamma_{j} \ln|q_{i} - q_{j}|^{2}.
\end{align*}
The vector field $X_{H}$ on $\R^{2N} \cong \C^{N}$ defined by the Hamiltonian system $\ins{X_{H}}\Omega = \d{H}$ then yields \eqref{eq:vortices}.

We note in passing that, strictly speaking, the system is not defined in the entire $\R^{2N} \cong \C^{N}$ because the Hamiltonian $H$ is not defined at those collision points, i.e., those with $\mathbf{x}_{i} = \mathbf{x}_{j}$ or $q_{i} = q_{j}$ with $i \neq j$.
However, we shall not be concerned with this issue here because this paper mainly concerns the stability of relative equilibria with positive inter-vortex distances.

It is straightforward to see that the system has $\SE(2) = \SO(2) \ltimes \R^{2}$-symmetry under the action
\begin{equation*}
  \label{eq:SE2-action}
  \SE(2) \times \C^{N} \to \C^{N};
  \qquad
  \bigl((e^{\rmi\theta}, a\bigr), \mathbf{q}) \mapsto e^{\rmi\theta}\mathbf{q} + a\mathbf{1},
\end{equation*}
where we identified $\R^{2}$ with $\C$, and defined $\mathbf{q} \defeq (q_{1}, \dots, q_{N})$ and $\mathbf{1} \defeq (1, \dots, 1)$, both in $\C^{N}$.
A \textit{relative equilibrium} of \eqref{eq:vortices} is a solution to \eqref{eq:vortices} whose trajectory is in an orbit of the above $\SE(2)$-action passing through the initial point $\mathbf{q}(0)$.
In other words, roughly speaking, a relative equilibrium is a solution in which the ``shape''---the relative positions and orientations of the $N$ vortices---does not change in time, although the individual vortices themselves may be moving in time.

We shall analyze the nonlinear stability of such relative equilibria by recasting relative equilibria into fixed points of the reduced equations obtained by $\SE(2)$-reduction of \eqref{eq:vortices}, specifically using the Lie--Poisson reduced dynamics developed in our previous work~\cite{Oh2019d}.

\subsection{Non-Equivariance and Stability of Relative Equilibria}
There are some known results on the stability of relative equilibria for Hamiltonian systems such as \citet{Mo1997} and \citet{PaRoWu2004} that do not require an explicit formulation of the reduced dynamics.

However, these results do not apply to our setting directly (see the next paragraph for extensions of these results to accommodate our case).
In fact, \citet{Mo1997} assumes that the symmetry group is compact, which is clearly not the case with $\SE(2)$ here.
\citet{PaRoWu2004} proved a non-compact extension of this result, but assume that the associated momentum map is equivariant.
That is not the case with our problem either in general: Consider the translational part of the $\SE(2)$-action \eqref{eq:SE2-action}, i.e., $\C \cong \R^{2}$-action on $\C^{N}$ as follows:
\begin{equation*}
  \C \times \C^{N} \to \C^{N};
  \qquad
  (a, \mathbf{q} \defeq (q_{1}, \dots , q_{N})) \mapsto \mathbf{q} + a\mathbf{1}.
\end{equation*}
where $\mathbf{1} \defeq (1, \dots, 1) \in \C^{N}$.
The momentum map associated with the translational subgroup $\R^{2}$-action of \eqref{eq:SE2-action} is essentially what is often called the linear impulse $\mathcal{I}(\mathbf{q}) \defeq -\rmi \sum_{j=1}^{N} \Gamma_{j} q_{j}$, and one easily sees that this momentum map is \textit{not} equivariant in general:
\begin{equation*}
  \mathcal{I}(\mathbf{q} + a\mathbf{1}) = \mathcal{I}(\mathbf{q}) - \rmi \Gamma a,
\end{equation*}
where $\Gamma$ is the total circulation:
\begin{equation*}
  \Gamma \defeq \sum_{i=1}^{N} \Gamma_{i}.
\end{equation*}
Specifically, one sees that the linear impulse $\mathcal{I}$ is equivariant if and only if $\Gamma = 0$.
However, the case with $\Gamma = 0$ is a very restrictive special case.
For example, there are many relative equilibria in the four-vortex problem with two pairs of vortices of equal circulations in which one pair has varying circulations and with $\Gamma \neq 0$; see \citet{HaRoSa2014}.
Some of those four-vortex relative equilibria are known to be linearly stable, but their nonlinear stability is an open question~\cite{Ro2013}.

We note in passing, however, that one may consider non-compact and non-equivariant extensions of the results of \citet{Mo1997,PaRoWu2004,LeSi1998} by, e.g., taking a special case assuming a free action in \citet{MoRo2011} (which implies $\mathfrak{h} = 0$ using the notation from \cite{MoRo2011}), considering a non-equivariant extension of \cite{LeSi1998,OrRa1999}, or considering a modified equivariant $\SE(2)$-action of \citet{Mo2000} (see also \citet{So1997})~\footnote{I would like thank one of the reviewers for bringing these ideas to my attention.}.

We also note that the stability of relative equilibria in this paper is restricted to the stability under perturbations that preserve the momentum map associated with the $\SE(2)$-action.
We refer to \cite{Pa1992,LeSi1998,OrRa1999} for stability of relative equilibria of Hamiltonian systems under more general perturbations again assuming the equivariance.
In particular, the result of \cite{Pa1992} implies that a relative equilibrium that is stable in our sense can still have drifts in the direction of the isotropy group of the initial value of the momentum map.

\subsection{Hamiltonian Formulation of Relative Dynamics}
\label{ssec:relative_dynamics-N=3}
In this paper, we shall take a different approach of performing stability analysis by using our previous work~\cite{Oh2019d} that gave an explicit Lie--Poisson formulation of the reduced dynamics.

Geometrically, the relative dynamics of the $N$-vortex system \eqref{eq:vortices} is the dynamics obtained by $\SE(2)$-symmetry reduction.
However, such a reduction is not particularly straightforward.
One would consider applying the semidirect product reduction (see, e.g., \cite[Sections~4.2]{MaMiOrPeRa2007}) by first applying the $\R^{2}$-reduction, followed by the $\SO(2)$-reduction.
In the $\R^{2}$-reduction, the non-equivariance of the linear impulse $\mathcal{I}$ mentioned above necessitates a separate handling of the cases with $\Gamma = 0$ and $\Gamma \neq 0$, but the resulting reduced space is rather simple~\cite{Oh2019d}.
For the $\SO(2)$-reduction, we used a dual pair to justify the formulation originally due to \citet{BoPa1998} and \citet{BoBoMa1999} (see also \cite{BoMa2015}) as the symplectic reduction.
The resulting reduced dynamics is formulated as a Lie--Poisson equation.
See also \cite{Du2013,ArDuNg2015,McMoVe2016b,DuMo2025} for the same idea applied to the $N$-body problem of celestial mechanics, and \cite{AlDu2020,DuSc2020} for the 3-body problem in $\R^{4}$.

For illustrative purposes, we shall first explain the reduced equations obtained in \cite{Oh2019d} (see also \cite{BoPa1998,BoBoMa1999}) for the special case with $N = 3$ and non-zero total circulation, i.e., $\Gamma \neq 0$; recall that the momentum map $\mathcal{I}$ is non-equivariant then.

The translational symmetry (the $\R^{2}$ subgroup symmetry of the above $\SE(2)$-symmetry) implies that the linear impulse $\mathcal{I} \defeq \sum_{j=1}^{3} \Gamma_{j} q_{j}$ (see, e.g., \citet[Section~2.1]{Ne2001} and \citet{Ar2007}) is an invariant of \eqref{eq:vortices}.
Hence we have
\begin{equation*}
  q_{3} = \frac{1}{\Gamma}(\mathcal{I} - \Gamma_{1}z_{1} - \Gamma_{2}z_{2})
  \quad\text{with}\quad
  z =
  \begin{bmatrix}
    z_{1} \\
    z_{2}
  \end{bmatrix}
  \defeq
  \begin{bmatrix}
    q_{1} - q_{3} \\
    q_{2} - q_{3}
  \end{bmatrix},
\end{equation*}
and so the dynamics essentially reduces to that of relative positions $z \in \C^{2}$ of vortices 1 and 2 with respect to vortex 3.
Now consider the matrix
\begin{equation}
  \label{eq:mu_N=3}
  \mu =
  \rmi
  \begin{bmatrix}
    \mu_{1} & \mu_{3} + \rmi \mu_{4} \\
    \mu_{3} - \rmi \mu_{4} & \mu_{2}
  \end{bmatrix}
  \defeq \rmi z z^{*} = \rmi
  \begin{bmatrix}
    |z_{1}|^{2} & z_{1} z_{2}^{*} \smallskip\\
    z_{2} z_{1}^{*} &  |z_{2}|^{2}
  \end{bmatrix}.
\end{equation}
We note in passing that this is essentially Eq.~(25) of \citet{Mo2015} used in the relative dynamics of the three-body problem of celestial mechanics.

\begin{figure}[hbtp]
  \centering
  \includegraphics[width=.45\linewidth]{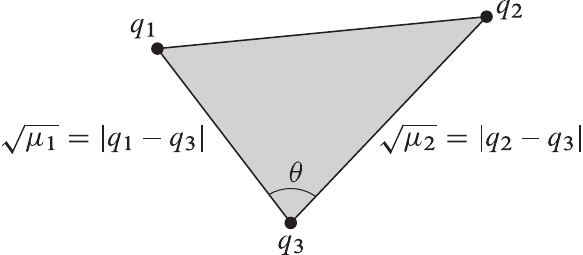}
  \caption{
    The entries $(\mu_{1}, \mu_{2}, \mu_{3})$ in the matrix $\mu$ in \eqref{eq:mu_N=3} give the variables that determine the triangular shape made by three vortices; see \eqref{eq:shape_variables-N=3}.
  }
  \label{fig:ShapeVariables}
\end{figure}

In what follows, we shall identify the matrix $\mu$ with vector $(\mu_{1}, \mu_{2}, \mu_{3}, \mu_{4}) \in \R^{4}$ whenever convenient.
We see that the three entries
\begin{equation}
  \label{eq:shape_variables-N=3}
  \begin{array}{c}
    \DS \mu_{1} \defeq |q_{1} - q_{3}|^{2},
    \qquad
    \mu_{2} \defeq |q_{2} - q_{3}|^{2},
    \medskip\\
    \mu_{3} \defeq \Re\bigl( (q_{1} - q_{3})(q_{2} - q_{3})^{*} \bigr)
    = |q_{1} - q_{3}| |q_{2} - q_{3}| \cos\theta
  \end{array}
\end{equation}
of $\mu$ determine the ``shape''---the triangle in this case---made by the three vortices as shown in \Cref{fig:ShapeVariables}.
So $\mu_{4}$ is a redundant variable to describe the shape.
Indeed, given the way $\mu$ is defined in terms of $z$ in \eqref{eq:mu_N=3}, one sees that $\rank \mu = 1$, and hence we have the constraint $\det \mu = \mu_{1} \mu_{2} - \mu_{3}^{2} - \mu_{4}^{2} = 0$.

It turns out that the time evolution of $\mu$ is governed by
\begin{equation}
  \label{eq:Lie-Poisson_N=3}
  \dot{\mu} = X_{h}(\mu) \defeq -\mu \frac{\delta h}{\delta\mu} \mathcal{K}^{-1} + \mathcal{K}^{-1} \frac{\delta h}{\delta\mu} \mu,
\end{equation}
where we defined
\begin{equation*}
  h(\mu) \defeq -\frac{1}{4\pi} \parentheses{
    \Gamma_{1} \Gamma_{3} \ln \mu_{1}
    + \Gamma_{2} \Gamma_{3} \ln \mu_{2}
    + \Gamma_{1} \Gamma_{2} \ln(\mu_{1} + \mu_{2} - 2\mu_{3})
  },
\end{equation*}
and
\begin{equation*}
  \fd{h}{\mu} \defeq \rmi
  \begin{bmatrix}
    2 \tpd{h}{\mu_{1}} & \tpd{h}{\mu_{3}} + \rmi\,\tpd{h}{\mu_{4}} \\
      \tpd{h}{\mu_{3}} - \rmi\,\tpd{h}{\mu_{4}} & 2 \tpd{h}{\mu_{2}}
  \end{bmatrix}
\end{equation*}
as well as
\begin{equation*}
  \mathcal{K} \defeq \frac{1}{\Gamma}
  \begin{bmatrix}
    -\Gamma_{1}(\Gamma - \Gamma_{1}) & \Gamma_{1} \Gamma_{2} \\
    \Gamma_{1} \Gamma_{2} & -\Gamma_{2}(\Gamma - \Gamma_{2})
  \end{bmatrix}.
\end{equation*}
The above equation is an instance of the so-called Lie--Poisson equation---a class of Hamiltonian systems defined on the dual of a Lie algebra; see, e.g., \citet[Chapter~13]{MaRa1999}, where $h$ is the Hamiltonian of the system---defined so that $h(\mu) = H(z)$ with $\mu = \rmi z z^{*}$.

Notice that \eqref{eq:Lie-Poisson_N=3} is defined on a vector space $\R^{4}$ as opposed to the shape sphere originally developed for the three-body problem in celestial mechanics by \citet{Mo2015} (see also \citet{Ro2018} for its application to the planar vortex problem).
While the increase in the dimension is certainly a drawback, it has the advantage of describing the dynamics explicitly using (global) coordinates on a vector space.

It is well known for $N = 3$ that if the three vortices are at the vertices of an equilateral triangle then it gives a relative equilibrium, i.e., the vortices move on the plane maintaining the ``shape'' of the initial equilateral triangle.
Without loss of generality, we may set the length of the edges to be 1 and assume that the vortices $(1,2,3)$ appear in this order counter-clockwise (the opposite of what is shown in \Cref{fig:ShapeVariables}).
Then, the corresponding values for the entries of $\mu$ are $\mu_{0} \defeq (1, 1, 1/2, -\sqrt{3}/2)$, written as a vector in $\R^{4}$.
It is a straightforward calculation to show that $X_{h}(\mu_{0}) = 0$, and hence $\mu_{0}$ is a fixed point of the above Lie--Poisson equation~\eqref{eq:Lie-Poisson_N=3}.
So one may analyze the stability of the equilateral triangle \textit{relative equilibria} by analyzing the stability of the corresponding \textit{fixed point} $\mu_{0}$ in the Lie--Poisson equation.
In \Cref{sec:Illustrative_Example}, we shall use this example again as an illustration of how our method works.

\subsection{Main Results and Outline}
The main result, \Cref{thm:Energy-Casimir}, gives a sufficient condition for Lyapunov stability of fixed points in the Lie--Poisson relative dynamics.
The result implies the stability of the corresponding relative equilibria in the planar $N$-vortex problem~\eqref{eq:vortices} to those perturbations that preserve the momentum map.
It provides a method to analyze nonlinear stability of relative equilibria, and is complementary to the existing methods~\cite{BaHo2016,MeRo2018,PeSaTa2015,Ro2013,Ro2018} that are mainly based on linear analysis, sometimes combined with Morse-theoretic methods.

The paper is organized to build up to the main result as follows:
In \Cref{sec:Illustrative_Example}, we first illustrate this approach using the equilateral triangle example from above, and recover the well-known stability condition by \citet{Sy1949} and \citet{Ar1979}.
In the rest of the paper, we shall build on this example and explain how the method used in this example may be generalized to the $N$-vortex case.

In \Cref{sec:Lie-Poisson}, we briefly review the Lie--Poisson relative dynamics for the planar $N$-vortex problem following our previous work~\cite{Oh2019d}, and then in \Cref{sec:Casimirs_and_Constraints}, we show that the Lie--Poisson relative dynamics possesses some invariants, and that it is constrained to some invariant submanifold.
These invariants and constraints are essential in overcoming the main drawback of our Lie--Poisson approach that it ostensibly increases the dimension of the dynamics, and play a critical role in the stability theorem, \Cref{thm:Energy-Casimir}.

In \Cref{sec:main_result}, we prove the main result, \Cref{thm:Energy-Casimir}.
It gives an Energy--Casimir-type stability condition that can be used explicitly for given relative equilibria.
Then, in \Cref{sec:examples}, we apply it to the relative equilibria of (i)~equilateral triangle with center with non-zero total circulation with varying circulation $\gamma$ for the center vortex, and (ii)~rhombus equilibria with two vortices of circulation 1 and the other two of circulation $\gamma \in (-1,1) \backslash \{0\}$.
\Cref{thm:Energy-Casimir} applied to the first case proves its nonlinear stability for $\gamma < -3$, in addition to the range of stability $-1/2 < \gamma < 1$ obtained by \citet{CaSc2000}.
In the latter case, we prove its nonlinear stability for $-2 + \sqrt{3} < \gamma < 0$; it was known to be only linearly stable due to \citet{Ro2013}.

\section{Illustrative Example with $N = 3$}
\label{sec:Illustrative_Example}
\subsection{Equilateral Triangle Relative Equilibrium}
We shall continue from \Cref{ssec:relative_dynamics-N=3} and analyze the stability of the equilateral triangle.
Its stability condition is well known; see, e.g., \citet{Sy1949} and \citet{Ar1979}; see also \citet[Theorem~2.2.2]{Ne2001} and \citet{Ro2013}.

The goal of this section is to give a preview of the main ideas of the paper by illustrating how our approach via the Lie--Poisson equation works in reproducing these results.

\subsection{Linear Stability of Equilateral Triangle}
Recall that we are analyzing the stability of the fixed point $\mu_{0} \defeq (1, 1, 1/2, -\sqrt{3}/2)$ of the Lie--Poisson equation~\eqref{eq:Lie-Poisson_N=3}.
By recasting \eqref{eq:Lie-Poisson_N=3} into a system in $\R^{4}$ and linearizing it around the fixed point $\mu_{0}$, we obtain a linear system in $\R^{4}$ with the coefficient matrix
\begin{equation*}
  A \defeq DX_{h}(\mu_{0}) =
  \frac{1}{4\pi}
  \begin{bmatrix}
    -2 \sqrt{3} \Gamma_{2} & 0 & 4 \sqrt{3} \Gamma_{2} & 0 \\
    0 & 2 \sqrt{3} \Gamma_{1} & -4 \sqrt{3} \Gamma_{1} & 0 \\
    -\sqrt{3} (\Gamma_{2}+\Gamma_{3}) & \sqrt{3} (\Gamma_{1}+\Gamma_{3}) & 2 \sqrt{3} (\Gamma_{2}-\Gamma_{1}) & 0 \\
    \Gamma_{2}-\Gamma_{3} & \Gamma_{3}-\Gamma_{1} & 2 (\Gamma_{1}-\Gamma_{2}) & 0 \\
  \end{bmatrix}.
\end{equation*}
Its eigenvalues are
\begin{equation*}
  \braces{ 0, 0, \pm\frac{\sqrt{3}}{2\pi} \sqrt{ -(\Gamma_{1} \Gamma_{2} + \Gamma_{1} \Gamma_{3} + \Gamma_{2} \Gamma_{3}) } }.
\end{equation*}
Hence we see that the fixed point $\mu_{0}$ is linearly unstable if $\Gamma_{1} \Gamma_{2} + \Gamma_{1} \Gamma_{3} + \Gamma_{2} \Gamma_{3} < 0$.

\subsection{Nonlinear Stability of Equilateral Triangle}
\label{ssec:equilateral}
We shall use the Energy--Casimir method of \citet{Ae1992}.
This method is particularly useful for Hamiltonian systems with additional invariants besides the Hamiltonian function.
For $N = 3$ and $\Gamma \defeq \sum_{i=1}^{3} \Gamma_{i} \neq 0$, the Lie--Poisson equation~\eqref{eq:Lie-Poisson_N=3} possesses two additional invariants:
\begin{equation*}
  C_{1}(\mu) \defeq \frac{\Gamma_{2} (\Gamma_{1}+\Gamma_{3}) \mu_{1} + \Gamma_{1} (\Gamma_{2}+\Gamma_{3}) \mu_{2} -2 \Gamma_{1} \Gamma_{2} \mu_{3}}{\Gamma}
\end{equation*}
as well as
\begin{equation}
  \label{eq:C_2-N=3}
  C_{2}(\mu) \defeq \mu_{1} \mu_{2} - \mu_{3}^{2} - \mu_{4}^{2},
\end{equation}
which was mentioned above as the constraint $\det\mu = 0$ for $\mu$ being rank-one.

The Energy--Casimir Theorem~\cite{Ae1992} applied to this setting says the following:
Suppose $DC_{1}(\mu_{0})$ and $DC_{2}(\mu_{0})$ are independent.
Then the fixed point $\mu_{0}$ of \eqref{eq:Lie-Poisson_N=3} is Lyapunov stable if there exist constants $c_{0}, c_{1}, c_{2} \in \R$ with $c_{0} \neq 0$ such that the function $f \defeq 4\pi c_{0} h + c_{1} C_{1} + c_{2} C_{2}$ satisfies the following:
\begin{enumerate}[(i)]
\item $Df(\mu_{0}) = 0$, and
  \smallskip
\item the Hessian $\mathsf{H} \defeq D^{2}f(\mu_{0})$ is positive definite on the tangent space at $\mu_{0}$ of the level set
  \begin{equation*}
    M \defeq \setdef{ \mu \in \R^{4} }{ C_{j}(\mu) = C_{j}(\mu_{0})\ \forall j \in \{1,2\} }
    = C_{1}^{-1}(C_{1}(\mu_{0})) \cap C_{2}^{-1}(0),
  \end{equation*}
  i.e., for every $v \in{\R}^{4} \backslash\{0\}$ satisfying $DC_{j}(\mu_{0}) \cdot v = 0$ with every $j \in \{1, 2\}$, one has $v^{T} \mathsf{H} v > 0$.
\end{enumerate}

We claim that if $\Gamma_{1} \Gamma_{2} + \Gamma_{1} \Gamma_{3} + \Gamma_{2} \Gamma_{3} > 0$ then $\mu_{0}$ is stable.
First, straightforward calculations show that the first condition is satisfied if $c_{1} = c_{0} \Gamma$ and $c_{2} = 0$.
So it remains to show that the second condition is satisfied with some choice of $c_{0} \neq 0$.

If $\Gamma_{1} \neq \Gamma_{2}$, then the tangent space $T_{\mu_{0}}M$ of the level set is spanned by
\begin{align*}
  v_{1} &\defeq \bigl( \sqrt{3}\Gamma_{2} (\Gamma_{1}+\Gamma_{3}),-\sqrt{3}\Gamma_{1} (\Gamma_{2}+\Gamma_{3}),0,\Gamma_{3} (\Gamma_{1}-\Gamma_{2}) \bigr), \\
  v_{2} &\defeq \bigl( \Gamma_{2}(\Gamma_{1}-\Gamma_{3}), \Gamma_{1} (\Gamma_{3}-\Gamma_{2}), \Gamma_{3} (\Gamma_{1}-\Gamma_{2}), 0 \bigr).
\end{align*}
Defining a $2 \times 2$ matrix $\mathcal{H}$ by setting $\mathcal{H}_{ij} \defeq v_{i}^{T} \mathsf{H} v_{j}$, the leading principal minors of $\mathcal{H}$ are
\begin{align*}
  d_{1} &= 3 c_{0} \Gamma_{1} \Gamma_{2} \Gamma_{3} (\Gamma_{1}+\Gamma_{2}) (\Gamma_{1}+\Gamma_{3}) (\Gamma_{2}+\Gamma_{3}), \\
  d_{2} &= 12 c_{0}^{2} \Gamma_{1}^{2} \Gamma_{2}^{2} \Gamma_{3}^4 (\Gamma_{1}-\Gamma_{2})^{2} (\Gamma_{1} \Gamma_{2} + \Gamma_{1} \Gamma_{3} + \Gamma_{2} \Gamma_{3})
\end{align*}
We see that $d_{2} > 0$ regardless of the value of $c_{0} \neq 0$ due to our assumptions.
For $d_{1}$, notice that $(\Gamma_{1}+\Gamma_{2}) (\Gamma_{1}+\Gamma_{3}) (\Gamma_{2}+\Gamma_{3}) \neq 0$ because if, say, $\Gamma_{2} = -\Gamma_{1}$, then $\Gamma_{1} \Gamma_{2} + \Gamma_{1} \Gamma_{3} + \Gamma_{2} \Gamma_{3} = -\Gamma_{1}^{2} < 0$, contradicting our assumption; similar contradiction if $\Gamma_{3} = -\Gamma_{1}$ or $\Gamma_{3} = -\Gamma_{2}$ as well.
Hence one may always choose $c_{0}$ so that its sign is the opposite of $(\Gamma_{1}+\Gamma_{2}) (\Gamma_{1}+\Gamma_{3}) (\Gamma_{2}+\Gamma_{3})$ to have $d_{1} > 0$.

If $\Gamma_{1} = \Gamma_{2} \neq \Gamma_{3}$, then the tangent space is spanned by
\begin{equation*}
  v_{1} = \bigl( 2\sqrt{3}\Gamma_{1}, 0, \sqrt{3}(\Gamma_{1} + \Gamma_{3}), 0, \Gamma_{3} - \Gamma_{1}),
  \qquad
  v_{2} = \bigl( \Gamma_{1}-\Gamma_{3}, \Gamma_{3}-\Gamma_{1}, 0, 0 \bigr),
\end{equation*}
and
\begin{equation*}
  d_{1} = 12 c_{0} \Gamma_{1}^{2} \Gamma_{3} (\Gamma_{1} + \Gamma_{3}),
  \qquad
  d_{2} = 12 c_{0}^{2} \Gamma_{1}^{2} \Gamma_{3}^{2} (\Gamma_{1} - \Gamma_{3})^{2} (\Gamma_{1} \Gamma_{2} + \Gamma_{1} \Gamma_{3} + \Gamma_{2} \Gamma_{3}),
\end{equation*}
and a similar argument shows that one can choose $c_{0}$ so that both $d_{1}$ and $d_{2}$ are positive.

If $\Gamma_{1} = \Gamma_{2} = \Gamma_{3}$, we have $v_{1} = (1, 0, 1, 0)$ and $v_{2} = (-1, 1, 0, 0)$, and $d_{1} = 2c_{0}\Gamma_{1}^{2}$ and $d_{2} = 3c_{0}^{2}\Gamma_{1}^{4}$.
So one may choose any $c_{0} > 0$.

Therefore, the Energy--Casimir Theorem implies that $\mu_{0}$ is stable if $\Gamma_{1} \Gamma_{2} + \Gamma_{1} \Gamma_{3} + \Gamma_{2} \Gamma_{3} > 0$.

\section{Lie--Poisson Relative Dynamics}
\label{sec:Lie-Poisson}
\subsection{$N$ Point Vortices with $\Gamma \neq 0$}
We may generalize the Lie--Poisson relative dynamics for $N = 3$ (\Cref{ssec:relative_dynamics-N=3}) to the general $N$-vortex case with
\begin{equation*}
  \Gamma \defeq \sum_{i=1}^{N} \Gamma_{i} \neq 0  
\end{equation*}
as is done in our previous work~\cite{Oh2019d}; see also \cite{BoPa1998,BoBoMa1999}.
Specifically, we shall use the relative coordinates
\begin{equation}
  \label{eq:z}
  z =
  \begin{bmatrix}
    z_{1} \\
    \vdots \\
    z_{N-1}
  \end{bmatrix}
  \defeq
  \begin{bmatrix}
    q_{1} - q_{N} \\
    \vdots \\
    q_{N-1} - q_{N}
  \end{bmatrix}
  \in \C^{N-1},
\end{equation}
and define
\begin{equation*}
  \mathbf{J}\colon \C^{N-1} \to \vK^{*}; 
  \qquad
  z \mapsto \rmi z z^{*},
\end{equation*}
where $\vK^{*}$ is the dual of the Lie algebra
\begin{equation}
  \label{eq:vK}
  \vK \defeq \setdef{ \xi \in \C^{(N-1)\times(N-1)} }{ \xi^{*} = -\xi }
\end{equation}
equipped with the (non-standard) Lie bracket
\begin{equation}
  \label{eq:Lie_bracket-v_K}
  [\xi, \eta]_{\mathcal{K}} \defeq \xi \mathcal{K}^{-1} \eta - \eta \mathcal{K}^{-1} \xi
\end{equation}
with
\begin{equation}
  \label{eq:K}
  \mathcal{K} \defeq \frac{1}{\Gamma}
  \begin{bmatrix}
    -\Gamma_{1}(\Gamma - \Gamma_{1}) & \Gamma_{1} \Gamma_{2} & \dots & \Gamma_{1} \Gamma_{N-1} \\
    \Gamma_{2} \Gamma_{1} & -\Gamma_{2}(\Gamma - \Gamma_{2}) & \dots & \Gamma_{2} \Gamma_{N-1} \\
    \vdots & \vdots & \ddots & \vdots \\
    \Gamma_{N-1} \Gamma_{1} & \Gamma_{N-1} \Gamma_{2} & \dots & -\Gamma_{N-1}(\Gamma - \Gamma_{N-1})
  \end{bmatrix},
\end{equation}
which is known to be invertible; see \cite[Remark~2.5]{Oh2019d}.

Notice that $\vK$, as a vector space (but not as a Lie algebra), is isomorphic to the unitary algebra $\u(N-1)$.
Then we may identify the dual $\vK^{*}$ with $\vK$ via the inner product
\begin{equation*}
  \ip{\xi}{\eta} \defeq \frac{1}{2}\tr(\xi^{*}\eta)
  \quad
  \forall \xi, \eta \in \vK,
\end{equation*}
and we shall use this identification in what follows.
Thus we may write an arbitrary element $\mu \in \vK^{*} \cong \vK$ as follows:
\begin{equation}
  \label{eq:mu}
  \mu = \rmi\,
  \begin{bNiceMatrix}
    \mu_{1} & \mu_{12} & \Cdots & \mu_{1,N-1} \\
    \mu_{12}^{*} & \Ddots & \Ddots & \vdots \\
    \vdots  & \Ddots & \Ddots  & \mu_{N-2,N-1} \\
    \mu_{1,N-1}^{*} & \Cdots & \mu_{N-2,N-1}^{*} & \mu_{N-1}
  \end{bNiceMatrix},
\end{equation}
where
\begin{equation*}
  \mu_{i} \in \R \text{ for } 1 \le i \le N-1,
  \qquad
  \mu_{ij} \in \C \text{ for } 1 \le i < j \le N-1.
\end{equation*}

Setting $\mu = \mathbf{J}(z)$, we see that
\begin{equation*}
  \mu_{i} = |z_{i}|^{2},
  \qquad
  \mu_{ij} = z_{i} z_{j}^{*},
\end{equation*}
and so the entries of $\mu$ are related to those variables that describe the relative configuration of the vortices or the ``shape'' made by the vortices.
What one can show is that $\mathbf{J}$ is a Poisson map with respect to the Poisson bracket associated with the symplectic form~\eqref{eq:Omega} and the Lie--Poisson bracket
\begin{equation*}
  \PB{f}{h}(\mu) \defeq \ip{\mu}{ \brackets{ \frac{\delta f}{\delta\mu}, \frac{\delta h}{\delta\mu} }_{\mathcal{K}} }
\end{equation*}
on $\vK$ using the Lie bracket~\eqref{eq:Lie_bracket-v_K}, where the derivative ${\delta f}/{\delta\mu} \in \vK$ is defined so that, for all $\mu,\nu \in \vK^{*}$,
\begin{equation*}
  \ip{\nu}{\frac{\delta f}{\delta\mu}}
  = \frac{1}{2}\tr\parentheses{ \nu^{*} \frac{\delta f}{\delta\mu} }
  = \left.\od{}{s}\right|_{s=0} f(\mu + s\nu).
\end{equation*}

Now, let us define the Hamiltonian $h\colon \mathfrak{v}_{K}^{*} \to \R$ so that $h \circ \mathbf{J} = H$ as follows:
\begin{equation}
  \label{eq:h}
  h(\mu) \defeq -\frac{1}{4\pi}\Biggl(
  \Gamma_{N} \sum_{i=1}^{N-1} \Gamma_{i} \ln \mu_{i}
  + \sum_{1\le i < j \le N} \Gamma_{i} \Gamma_{j} \ln(\mu_{i} + \mu_{j} - 2\Re\mu_{ij})
  \Biggr).
\end{equation}
\begin{subequations}
  Then the dynamics of $\mu$ is governed by the Lie--Poisson equation
  \label{IVP:Lie-Poisson}
  \begin{equation}
  \label{eq:Lie-Poisson}
    \dot{\mu} = X_{h}(\mu) \defeq -\ad_{\delta h/\delta\mu}^{*}\mu \defeq -\mu \frac{\delta h}{\delta\mu} \mathcal{K}^{-1} + \mathcal{K}^{-1} \frac{\delta h}{\delta\mu} \mu,
  \end{equation}
  and the initial condition
  \begin{equation}
    \label{IC:Lie-Poisson}
    \mu(0) = \mathbf{J}(z(0))
    \quad\text{with}\quad
    z(0) = \begin{bmatrix}
      q_{1}(0) - q_{N}(0) \\
      \vdots \\
      q_{N-1}(0) - q_{N}(0)
    \end{bmatrix}.
  \end{equation}
\end{subequations}
It is proved in \cite[Corollary~3.6]{Oh2019d} that the above initial value problem~\eqref{IVP:Lie-Poisson} indeed gives the $\SE(2)$-reduced dynamics of \eqref{eq:vortices} obtained by symplectic reduction, assuming non-zero angular impulse.

\section{Casimir Invariants and Constraints}
\label{sec:Casimirs_and_Constraints}
\subsection{Casimirs and Constraints}
The advantage of the above Lie--Poisson formulation of the relative dynamics is that the reduced dynamics is defined on a vector space $\vK$.
As mentioned earlier, this is not the case if one performs the $\SE(2)$-reduction using the standard techniques of symplectic/Hamiltonian reduction.

However, there is one obvious drawback as well:
The (real) dimension of $\vK$ is $(N - 1)^{2}$, and is far greater than what one would expect for the reduced dynamics, especially for large $N$.
This is the price we pay for reformulating the reduced dynamics by embedding it into a vector space.

Nevertheless, the increase in the dimension may be mitigated by the presence of invariants and constraints in our formulation.
As we shall see, this is particularly helpful for our purpose because these constraints and invariants, used in the Energy--Casimir-type method as in \Cref{ssec:equilateral}, effectively reduce the dimension of the problem.
In fact, recall that, for the $N = 3$ case from \Cref{ssec:equilateral}, $\mu$ had four independent entries including $\mu_{4}$, which was redundant for describing the triangular shape.
However, with the help of the Energy--Casimir Theorem along with invariants $C_{1}$ and $C_{2}$, the stability analysis was reduced to the definiteness of a $2 \times 2$ matrix.

Let us describe the constraints for the relative dynamics~\eqref{IVP:Lie-Poisson} alluded above.
Again we shall assume that $\Gamma \neq 0$ because $\Gamma = 0$ is the case where the momentum map is equivariant, and may be treated using a result from \citet{PaRoWu2004}.

First recall that an arbitrary element $\mu \in \vK$ is written as an $(N-1) \times (N-1)$ skew-Hermitian matrix as shown in \eqref{eq:vK} and \eqref{eq:mu}.
However, recall also that \eqref{eq:Lie-Poisson} is obtained via $\mu = \mathbf{J}(z) = \rmi z z^{*}$, where $z$ is the positions of vortices $1$ to $N-1$ relative to vortex $N$ as defined in \eqref{eq:z}.
Therefore, one sees that $\rank \mu = 1$, i.e., $\mu$ is in the following subset of $\vK$:
\begin{align*}
  \vK^{(1)} &\defeq \setdef{ \mu \in \vK }{ \rank \mu = 1 } \\
  &= \setdef{ \mu \in \C^{(N-1)\times(N-1)} }{ \mu^{*} = -\mu,\ \rank \mu = 1 }.
\end{align*}
Moreover, notice that we also have
\begin{equation*}
  \mu_{i} = |z_{i}|^{2},
  \qquad
  \mu_{ij} = z_{i} z_{j}^{*},
\end{equation*}
and so, assuming no collisions of the vortices, every entry of the matrix $\mu$ is non-zero.
Hence $\mu$ is also in the following open subset of $\vK$:
\begin{equation}
  \label{eq:nonzerovK}
  \nonzerovK
  \defeq \setdef{ \mu \in \vK }{ \mu_{i} \neq 0 \text{ for } 1 \le i \le N-1 \text{ and } \mu_{ij} \neq 0 \text{ for } 1 \le i < j \le N-1 }.
\end{equation}
Therefore, the matrix $\mu$ in the relative dynamics~\eqref{eq:Lie-Poisson} lives in the following subset of $\vK$:
\begin{equation}
  \label{eq:nonzerovK1}
  \nonzerovK^{(1)} \defeq \vK^{(1)} \cap \nonzerovK.
\end{equation}

We shall summarize the invariants and constraints of the relative dynamics~\eqref{IVP:Lie-Poisson} as follows:
\begin{proposition}[Casimirs and constraints of Lie--Poisson relative dynamics]
  \label[proposition]{prop:Casimirs_and_invariants}
  The Lie--Poisson relative dynamics~\eqref{eq:Lie-Poisson} of $N$ point vortices with non-varnishing total circulation ($\Gamma \neq 0$) possesses Casimirs (invariants)
  \begin{equation}
    \label{eq:Casimirs}
    C_{j}(\mu) \defeq \tr((\rmi\,\mathcal{K} \mu)^{j})
    \quad
    j \in \N,
  \end{equation}
  where $C_{j}$ with $j \ge N$ are dependent on $\{ C_{i} \}_{i=1}^{N-1}$.
  Moreover, assuming no collisions of the vortices, the subset $\nonzerovK^{(1)}$ defined in \eqref{eq:nonzerovK1} is an invariant set of the relative dynamics governed by \eqref{eq:Lie-Poisson}.
\end{proposition}
\begin{proof}
  That $\{ C_{j} \}_{j \in \N}$ are Casimirs (invariants) is proved in our previous work~\cite[Corollary~3.6]{Oh2019d}.
  One can prove the dependence in the same way as in the proof of \cite[Proposition~3.6]{Oh2023a} for the relative dynamics of vortices on the sphere.
  
  So it remains to show that $\nonzerovK^{(1)}$ is an invariant set.
  Indeed, let $t \mapsto q(t)$ be the solution of the initial value problem of the original evolution equations~\eqref{eq:vortices-C}.
  Letting $t \mapsto z(t)$ be the corresponding $z$ via \eqref{eq:z}, we have $z_{i}(t) \defeq q_{i}(t) - q_{N}(t) \neq 0$ for every $i \in \{1, \dots, N-1\}$ and every $t$ assuming no collisions.
  Then, setting $\mu(0) = \rmi\,z(0) z(0)^{*}$, both $t \mapsto \mu(t)$ and $t \mapsto \mathbf{J}(z(t)) = \rmi\,z(t) z(t)^{*}$ satisfy the initial value problem~\eqref{eq:Lie-Poisson}.
  Hence by uniqueness we have $\mu(t) = \rmi\,z(t) z(t)^{*}$, and so $\rank \mu(t) = 1$, i.e., $\mu(t) \in \vK^{(1)}$ for every $t$.
  However, since $z_{i}(t) \neq 0$ for every $i \in \{1, \dots, N-1\}$ and every $t$, we see that
  \begin{align*}
    \mu_{i}(t) &= |z_{i}(t)|^{2} \neq 0 \text{ for } 1 \le i \le N-1, \\
    |\mu_{ij}(t)| &= |z_{i}(t)| |z_{j}(t)| \neq 0  \text{ for } 1 \le i < j \le N - 1,
  \end{align*}
  that is $\mu(t) \in \nonzerovK$ for every $t$ as well.
\end{proof}

\subsection{Imposing Constraints}
For stability analysis, it is desirable to impose the constraint to $\nonzerovK^{(1)}$ in a more explicit manner.
To that end, let us first prove some basic lemma on such matrices:
\begin{lemma}
  Let $n \ge 2$ and $A \in \C^{n\times n}$ be a Hermitian matrix with no vanishing elements.
  Then $\rank A = 1$ if and only if all the determinants of the $2 \times 2$ submatrices sweeping the upper triangular part and the subdiagonal (see the picture below) vanish.
  \begin{equation*}
    A = 
    \begin{bmatrix}
      \tikzmarknode{11}{a_{11}} & \tikzmarknode{12}{a_{12}} & \tikzmarknode{13}{a_{13}} & \tikzmarknode{14}{\vphantom{a_{14}}\cdots} & \tikzmarknode{15}{\vphantom{a_{15}}\cdots} & \tikzmarknode{16}{a_{1n}}  \\ 
      \tikzmarknode{21}{a_{21}} & \tikzmarknode{22}{a_{22}} & \tikzmarknode{23}{a_{23}} & \tikzmarknode{24}{\vphantom{a_{24}}\cdots} & \tikzmarknode{25}{\vphantom{a_{25}}\cdots} & \tikzmarknode{26}{a_{2n}} \\
      \tikzmarknode{31}{a_{31}} & \tikzmarknode{32}{a_{32}} & \tikzmarknode{33}{a_{33}} & \tikzmarknode{34}{\vphantom{a_{34}}\cdots} & \tikzmarknode{35}{\vphantom{a_{35}}\cdots} & \tikzmarknode{36}{a_{3n}} \\
      \tikzmarknode{41}{a_{41}} & \tikzmarknode{42}{a_{42}} & \tikzmarknode{43}{a_{43}} & \tikzmarknode{44}{\vphantom{a_{44}}\raisebox{0pt}[0pt][0pt]{$\ddots$}} & \tikzmarknode{45}{\vphantom{a_{45}}\cdots} & \tikzmarknode{46}{a_{4n}} \\
      \vdots & \vdots & \vdots & \vdots & \tikzmarknode{55}{\vphantom{\raisebox{3pt}{$a_{55}$}}\raisebox{0pt}[0pt][0pt]{$\ddots$}} & \tikzmarknode{56}{\vphantom{\raisebox{3pt}{$a_{56}$}}\raisebox{0pt}[0pt][0pt]{$\vdots$}} \\
      a_{n1} & a_{n2} & a_{n3} & \dots & \cdots & \tikzmarknode{66}{a_{nn}}  \\
    \end{bmatrix}
  \end{equation*}
  \AddToShipoutPictureBG*{%
    \begin{tikzpicture}[overlay,remember picture]
      \node[fit=(11) (22),draw=red,fill=red,fill opacity=0.5]{};
      \node[fit=(12) (23),draw=red,fill=red,fill opacity=0.5]{};
      \node[fit=(13) (24),draw=red,fill=red,fill opacity=0.5]{};
      \node[fit=(14) (25),draw=red,fill=red,fill opacity=0.5]{};
      \node[fit=(15) (26),draw=red,fill=red,fill opacity=0.5]{};
      \node[fit=(22) (33),draw=red,fill=red,fill opacity=0.5]{};
      \node[fit=(23) (34),draw=red,fill=red,fill opacity=0.5]{};
      \node[fit=(24) (35),draw=red,fill=red,fill opacity=0.5]{};
      \node[fit=(25) (36),draw=red,fill=red,fill opacity=0.5]{};
      \node[fit=(33) (44),draw=red,fill=red,fill opacity=0.5]{};
      \node[fit=(34) (45),draw=red,fill=red,fill opacity=0.5]{};
      \node[fit=(35) (46),draw=red,fill=red,fill opacity=0.5]{};
      \node[fit=(44) (55),draw=red,fill=red,fill opacity=0.5]{};
      \node[fit=(45) (46) (56),draw=red,fill=red,fill opacity=0.5]{};
      \node[fit=(55) (66),draw=red,fill=red,fill opacity=0.5]{};
    \end{tikzpicture}%
  }
\end{lemma}
\begin{proof}
  The sufficiency is straightforward, because $\rank A = 1$ implies that all the rows of $A$ are proportional to each other, and hence all the determinants in question vanish.

  Let us prove the necessity by induction on $n$.
  It clearly holds for $n = 2$ because $A \neq 0$ by assumption: $\rank A = 1$ if and only if $\det A = 0$.
  Suppose that the assertion holds for $n - 1$, and we shall show that it holds for $n$.
  To that end, suppose that all the determinants indicated in the statement for $A \in \C^{n\times n}$ vanish.
  Deleting the last column and the last row of $A$ to have $A_{n-1} \in \C^{(n-1)\times(n-1)}$.
  Then those $2 \times 2$ determinants sweeping the upper triangular part and the subdiagonal of $A_{n-1}$ vanish by assumption.
  Thus the induction hypothesis implies that $\rank A_{n-1} = 1$, and hence all the rows of $A_{n-1}$ are proportional to each other.
  The remaining $n-1$ vanishing determinants, sweeping the last two columns of $A$, impose that the last two columns of $A$ are proportional to each other, because all elements of $A$ are non-zero.
  As a result, all the rows of $A$ except the last one are proportional to each other.
  However, since $A$ is Hermitian, the last two columns of $A$ being proportional to each other implies that the last two rows of $A$ are proportional to each other as well.
  Therefore, all $n$ rows of $A$ are proportional to each other, and thus $\rank A = 1$ because $A \neq 0$.
\end{proof}

Let us apply the above lemma to the Hermitian matrix
\begin{equation}
  \label{eq:-rmi_mu}
  -\rmi\mu =
  \begin{bNiceMatrix}
    \mu_{1} & \mu_{12} & \Cdots & \mu_{1,N-1} \\
    \mu_{12}^{*} & \Ddots & \Ddots & \vdots \\
    \vdots  & \Ddots & \Ddots  & \mu_{N-2,N-1} \\
    \mu_{1,N-1}^{*} & \Cdots & \mu_{N-2,N-1}^{*} & \mu_{N-1}
  \end{bNiceMatrix}.
\end{equation}
Those determinants from the above lemma applied to it yield the following functions:
\begin{equation}
  \label{eq:Fs}
  \begin{split}
    F_{i}&\colon \nonzerovK \to \R;
           \quad
           F_{i}(\mu) \defeq
           \begin{vmatrix}
             \mu_{i} & \mu_{i,i+1} \\
             \mu_{i,i+1}^{*} & \mu_{i+1}
           \end{vmatrix}
           \quad\text{for}\quad
           1 \le i \le N-2,
    \\
    F_{ij}&\colon \nonzerovK \to \C;
            \quad
            F_{ij}(\mu) \defeq
            \begin{vmatrix}
              \mu_{i,j} & \mu_{i,j+1} \\
              \mu_{i+1,j} & \mu_{i+1,j+1}
            \end{vmatrix}
            \quad\text{for}\quad
            1 \le i < j \le N-2,
  \end{split}
\end{equation}
where we note that $F_{ij}$ with $j = i+1$ have diagonal elements $\mu_{i+1,i+1} = \mu_{i+1} \in \R$ in them:
\begin{equation*}
  F_{i,i+1} =
  \begin{vmatrix}
    \mu_{i,i+1} & \mu_{i,i+2} \\
    \mu_{i+1,i+1} & \mu_{i+1,i+2}
  \end{vmatrix}
  =
  \begin{vmatrix}
    \mu_{i,i+1} & \mu_{i,i+2} \\
    \mu_{i+1} & \mu_{i+1,i+2}
  \end{vmatrix}.
\end{equation*}
Combining all of the above functions to define
\begin{equation}
  \label{eq:F}
  \begin{split}
    F \colon &\nonzerovK \to \R^{N-1} \times \C^{N-1 \choose 2} \cong \R^{(N - 2)^{2}}; \\
             & F(\mu) \defeq (F_{1}(\mu), \dots, F_{N-2}(\mu), F_{12}(\mu), \dots, F_{N-3,N-2}(\mu)),
  \end{split}
\end{equation}
where we identified $\C$ with $\R^{2}$ in the standard manner, and note that $N - 2 + 2 \cdot {N - 2 \choose 2} = (N - 2)^{2}$.

Now we may characterize the subset $\nonzerovK^{(1)}$---to which $\mu$ is constrained---in a more effective manner using the above function $F$:
\begin{proposition}
  \label[proposition]{prop:constraint_functions}
  For every $N \ge 3$, the following hold:
  \begin{enumerate}[(i)]
  \item $\nonzerovK^{(1)} = F^{-1}(0)$.
    \smallskip
  \item $F$ is a submersion.
    \label{prop:constraint_functions-submersion}
    \smallskip
  \item $\nonzerovK^{(1)}$ is a $(2N - 3)$-dimensional submanifold of $\nonzerovK$.
  \end{enumerate}
\end{proposition}
\begin{proof}
  \leavevmode
  \begin{enumerate}[(i)]
  \item It follows immediately from the above lemma applied to $-\rmi \mu$ with $\mu \in \nonzerovK$.
    \smallskip
  \item See \Cref{sec:submersion} for a tedious but elementary proof.
    \smallskip
  \item Since $F$ is a submersion, a standard result from differential geometry implies that $F^{-1}(0)$ is a submanifold of codimension $(N - 2)^{2}$ of $\nonzerovK$.
    Notice from the definition \eqref{eq:nonzerovK} of $\nonzerovK$ that it is an open subset of $\vK$.
    Therefore, the (real) dimension of $\nonzerovK$ is given by
    \begin{equation*}
      \dim \nonzerovK = \dim \vK = (N - 1) + 2 \cdot \frac{(N - 1)^{2} - (N - 1)}{2}  = (N - 1)^{2},
    \end{equation*}
    because, for $\mu \in \vK$, the diagonal entries of $-\rmi\mu$ (see \eqref{eq:-rmi_mu}) are real, whereas its off-diagonal ones are complex; see \eqref{eq:vK} for the definition of $\vK$.
    Hence we have
    \begin{equation*}
      \dim F^{-1}(0) = (N - 1)^{2}- (N - 2)^{2} = 2N - 3. \qedhere
    \end{equation*}
  \end{enumerate}
\end{proof}

\begin{remark}
  \label{rem:Rs_are_not_invariants}
  The functions $\{ F_{i} \}_{i=1}^{N-2}$ and $\{ F_{ij} \}_{1\le i < j \le N-2}$ are in general \textit{not} invariants of the Lie--Poisson vector field $X_{h}$ from \eqref{eq:Lie-Poisson}, i.e., the Lie derivatives $X_{h}[F_{i}]$ and $X_{h}[F_{ij}]$ do \textit{not} vanish in general.
  However, the \textit{zero level set} of each of those functions is an \textit{invariant set} of the dynamics defined by $X_{h}$, as \Cref{prop:Casimirs_and_invariants,prop:constraint_functions} imply.
  The Casimir $C_{2}$ from $N = 3$ (\Cref{ssec:equilateral}) is nothing but $F_{1}$, but this is a special case in which $F_{1}$ becomes an invariant of the dynamics.
\end{remark}

\section{Stability Theorem for Fixed Points in Lie--Poisson Relative Dynamics}
\label{sec:main_result}
\subsection{Stability Lemma}
As proved in \Cref{prop:Casimirs_and_invariants,prop:constraint_functions}, the zero level set $\nonzerovK^{(1)} = F^{-1}(0)$ is an invariant set of the relative dynamics~\eqref{IVP:Lie-Poisson}.
However, as mentioned in \Cref{rem:Rs_are_not_invariants}, the functions $\{ F_{i} \}_{i=1}^{N-2}$ and $\{ F_{ij} \}_{1\le i < j \le N-2}$ defined in \eqref{eq:Fs} comprising $F$ are \textit{not} invariants of the Lie--Poisson relative dynamics in general.
Hence they are not on the same footing as the Casimirs $C_{1}$ and $C_{2}$ used in \Cref{ssec:equilateral} to prove the stability of the equilateral triangle.
This necessitates a generalization of the Energy--Casimir Theorem of \citet{Ae1992} that incorporates the invariant level set $F^{-1}(0)$.

To that end, let us first identify $\vK$ with $\R^{(N-1)^{2}}$ so that $\nonzerovK$ is an open subset of $\R^{(N-1)^{2}}$.
Then we have the following modification of a result from \citet[Theorem on p.~326]{Ae1992}:

\begin{lemma}[Stability lemma in invariant level set $F^{-1}(0)$]
  \label{lem:stability}
  Consider the Lie--Poisson relative dynamics~\eqref{eq:Lie-Poisson} on the invariant level set $\nonzerovK^{(1)} = F^{-1}(0)$ with Hamiltonian $C_{0} \defeq h$ and Casimirs $\{ C_{j} \}_{j=1}^{N-1}$, and let $\mu_{0} \in F^{-1}(0)$ be a fixed point of the dynamics.
  Suppose that there exists an open ball $\mathbb{B}_{r}(\mu_{0}) \subset \R^{(N-1)^{2}}$ centered at $\mu_{0}$ with radius $r > 0$ such that the level set of $\{ C_{j} \}_{j=0}^{N-1}$ at the value $(C_{0}(\mu_{0}), \dots, C_{N-1}(\mu_{0})) \in \R^{N}$ contains no points in $\mathbb{B}_{r}(\mu_{0}) \cap F^{-1}(0)$ other than $\mu_{0}$ itself, i.e.,
  \begin{equation}
    \label{eq:levelset}
    \setdef{ \mu \in \mathbb{B}_{r}(\mu_{0}) \cap F^{-1}(0) }{
      C_{j}(\mu) =  C_{j}(\mu_{0})
      \
      \forall
      j \in \{0, \dots, N-1\}
    }
    = \{ \mu_{0} \}.
  \end{equation}
  Then $\mu_{0}$ is Lyapunov stable to perturbations that preserve $F^{-1}(0)$.
\end{lemma}
\begin{proof}
  We essentially follow the proof of \cite[Theorem on p.~326]{Ae1992} with an adjustment to take care of the level set $F^{-1}(0)$.
  In what follows, we shall denote the flow of \eqref{eq:Lie-Poisson} by $\Phi$, i.e., for every $\nu \in F^{-1}(0)$, $t \mapsto \Phi_{t}(\nu)$ gives the solution to \eqref{eq:Lie-Poisson} with initial condition $\mu(0) = \nu$.
  
  By contradiction.
  Suppose that $\mu_{0} \in F^{-1}(0)$ is not Lyapunov stable.
  Then there exists an open neighborhood $U \subset F^{-1}(0)$ of $\mu_{0}$ with the following property:
  For every open neighborhood $V \subset U$ of $\mu_{0}$, there exists $\nu_{0} \in V$ such that the solution $t \mapsto \Phi_{t}(\nu_{0})$ eventually exits $U$, i.e., $\Phi_{\tau}(\nu_{0}) \notin U$ for some $\tau > 0$.

  Since $F^{-1}(0)$ is a submanifold of $\nonzerovK$, there exists an open subset $\mathcal{U}$ of $\nonzerovK$ such that $U = \mathcal{U} \cap F^{-1}(0)$.
  However, since $\nonzerovK$ in turn is an open subset of $\R^{(N-1)^{2}}$, there exists $\eps \in (0,r)$ such that $\mathbb{B}_{2\eps}(\mu_{0}) \subset \mathcal{U}$.
  Hence the closed ball $\bar{\mathbb{B}}_{\eps}(\mu_{0})$ is contained in $\mathcal{U}$, implying that $\bar{\mathbb{B}}_{\eps}(\mu_{0}) \cap F^{-1}(0) \subset U$.

  Then, for every $n \in \N$, we may take $\mathbb{B}_{\eps/n}(\mu_{0}) \cap F^{-1}(0)$ as the neighborhood $V$ mentioned above to find $\mu_{0}^{(n)} \in \mathbb{B}_{\eps/n}(\mu_{0}) \cap F^{-1}(0)$ such that $\Phi_{\tau}(\mu_{0}^{(n)}) \notin U$ for some $\tau > 0$.
  But then this implies that there exists $t_{n} > 0$ such that $\lambda_{n} \defeq \Phi_{t_{n}}(\mu_{0}^{(n)}) \in \mathbb{S}_{\eps}(\mu_{0}) \cap F^{-1}(0)$, where $\mathbb{S}_{\eps}(\mu_{0}) = \partial\bar{\mathbb{B}}_{\eps}(\mu_{0})$ is the sphere in $\R^{(N-1)^{2}}$ centered at $\mu_{0}$ with radius $\eps$.

  Clearly $\mathbb{S}_{\eps}(\mu_{0}) \cap F^{-1}(0)$ is compact, and so the sequence $\{ \lambda_{n} \}_{n=1}^{\infty}$---after passing down to a subsequence (without renumbering for simplicity)---converges to $\lambda \in \mathbb{S}_{\eps}(\mu_{0}) \cap F^{-1}(0)$.
  However, since $\{ C_{0} \defeq h \} \cup \{ C_{j} \}_{j=1}^{N-1}$ are invariants of the dynamics \eqref{eq:Lie-Poisson} and are continuous, we have
  \begin{equation*}
    \lim_{n\to\infty} C_{j}( \mu_{0}^{(n)} ) = \lim_{n\to\infty} C_{j}( \lambda_{n}) = C_{j}(\lambda)
    \quad
    \forall
    j \in \{0, \dots, N-1\}.
  \end{equation*}
  On the other hand, we also have $\mu_{0}^{(n)} \to \mu_{0}$ as $n \to \infty$, and so again by continuity,
  \begin{equation*}
    \lim_{n\to\infty} C_{j}( \mu_{0}^{(n)} ) = C_{j}(\mu_{0})
    \quad
    \forall
    j \in \{0, \dots, N-1\},
  \end{equation*}
  implying that
  \begin{equation*}
    C_{j}(\lambda) = C_{j}(\mu_{0})
    \quad
    \forall
    j \in \{0, \dots, N-1\}.
  \end{equation*}
  
  As a result, we see that $\mu_{0}$ and $\lambda$ are in the same level set, but they are distinct because $\lambda \in \mathbb{S}_{\eps}(\mu_{0})$.
  However, since $\eps < r$, we see that $\lambda \in \mathbb{B}_{r}(\mu_{0}) \cap F^{-1}(0)$ as well, contradicting the assumption~\eqref{eq:levelset} that $\mu_{0}$ is the only point in the level set.
\end{proof}

\subsection{Energy--Casimir Theorem for Lie--Poisson Relative Dynamics}
Now we shall build on the above stability lemma to prove an Energy--Casimir theorem in the invariant set $F^{-1}(0)$.
\begin{theorem}[Energy--Casimir method for Lie--Poisson relative dynamics]
  \label{thm:Energy-Casimir}
  Let $\mu_{0} \in F^{-1}(0) = \nonzerovK^{(1)}$ be a fixed point of the Lie--Poisson dynamics~\eqref{eq:Lie-Poisson}, and $\{ C_j \}_{j=1}^{K}$ be a subset of the Casimirs $\{ C_j \}_{j=1}^{N-1}$ such that
  \begin{equation*}
    \{ C_j \}_{j=1}^{K} \cup \{ F \} = \{ C_j \}_{j=1}^{K} \cup \{ F_{i} \}_{i=1}^{N-2} \cup \{ \Re F_{ij}, \Im F_{ij} \}_{1\le i < j \le N-2}
  \end{equation*}
  are independent at $\mu_{0}$ in the sense that their exterior differentials at $\mu_{0}$ are linearly independent.
  Suppose that there exist $a_{0} \in \R\backslash\{0\}$, $\{ a_{i} \in \R \}_{i=1}^{K}$ and $\{ b_{i} \in \R \}_{i=1}^{N-2} \cup \{ (c_{ij}, d_{ij}) \in \R^{2} \}_{1\le i < j \le N-2}$ such that
  \begin{align*}
    f(\mu) &\defeq a_{0} h(\mu) + \sum_{i=1}^{K} a_{i} C_{i}(\mu) + \sum_{i=1}^{N-1} b_{i} F_{i}(\mu) \\
    &\qquad + \sum_{1\le i < j \le N-2} (c_{ij} \Re F_{ij}(\mu) + d_{ij} \Im F_{ij}(\mu))
  \end{align*}
  satisfies the following:
  \begin{enumerate}[(i)]
  \item $Df(\mu_{0}) = 0$; and
    \label{cond:Energy-Casimir-1}
    \smallskip
  \item the Hessian $\mathsf{H} \defeq D^{2}f(\mu_{0})$ is positive definite on the tangent space at $\mu_{0}$ of the submanifold
    \begin{align*}
      M &\defeq \setdef{ \mu \in \nonzerovK }{ F(\mu) = 0,\ C_{j}(\mu) = C_{j}(\mu_{0})\ \forall j \in \{1, \dots, K\} } \\
        &= F^{-1}(0) \cap \parentheses{ \bigcap_{j=1}^{K} C_{j}^{-1}(C_{j}(\mu_{0})) }
    \end{align*}
    i.e., $v^{T} \mathsf{H} v > 0$ for every $v \in (T_{\mu_{0}}M)\backslash\{0\} \subset \R^{(N-1)^{2}}$.
    \label{cond:Energy-Casimir-2}
  \end{enumerate}
  Then $\mu_{0}$ is Lyapunov stable to perturbations that preserve $F^{-1}(0)$.
\end{theorem}
\begin{remark}
  One may instead have negative definiteness in condition~\eqref{cond:Energy-Casimir-2}, but it implies the positive definiteness by having $-f$ in place of $f$, i.e., by changing the signs of all the coefficients without breaking condition~\eqref{cond:Energy-Casimir-1}.
  We prefer to have positive-definiteness because it is easier to use Sylvester's criterion that way.
\end{remark}
\begin{remark}
  Stability of $\mu_{0} \in F^{-1}(0) = \nonzerovK^{(1)}$ to perturbations that preserve $F^{-1}(0)$ imply the stability of the relative equilibria to perturbations that preserve the momentum map.
  This is because the dual pair from \cite{Oh2019d} implies that each level set of the momentum map is mapped by $\mathbf{J}$ to a coadjoint orbit in $\mathfrak{v}_{\mathcal{K}}$, and the coadjoint orbit is contained in $F^{-1}(0)$ (assuming no collisions) as discussed in the proof of \Cref{prop:Casimirs_and_invariants}.
  Hence every momentum-map-preserving perturbation is $F^{-1}(0)$-preserving.
\end{remark}
\begin{proof}[Proof of \Cref{thm:Energy-Casimir}]
  It suffices to show that the above conditions are sufficient for $\mu_{0}$ to satisfy the condition stated in \Cref{lem:stability}.

  Recall first that $\nonzerovK$ is an open subset of $\vK \cong \R^{(N-1)^{2}}$.
  By the standard second-order condition in  equality-constrained optimization (see, e.g., \citet[Theorem~5.4 on p.~118]{Su1996}), assumptions~\eqref{cond:Energy-Casimir-1} and \eqref{cond:Energy-Casimir-2} imply that $a_{0}h$ achieves a strict local minimum at $\mu_{0}$ in $M$; hence $h$ itself achieves either strict local minimum or maximum at $\mu_{0}$ in $M$, because $a_{0} \neq 0$.
  But then this implies that there exists an open ball $\mathbb{B}_{r}(\mu_{0}) \subset \R^{(N-1)^{2}}$ such that the assumption \eqref{eq:levelset} in \Cref{lem:stability} holds.
  Therefore, by \Cref{lem:stability}, $\mu_{0}$ is Lyapunov stable.
\end{proof}

\section{Examples}
\label{sec:examples}
\subsection{Stability of Equilateral Triangle with Center}
\begin{figure}[htbp]
  \centering
  \includegraphics[width=.4\textwidth]{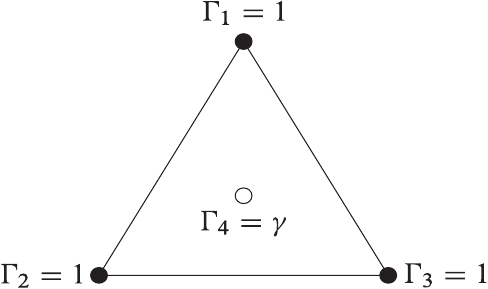}
  \caption{
    The relative equilibrium in the shape of equilateral triangle with center was shown to be Lyapunov stable if or $-1/2 < \gamma < 1$ in \cite{CaSc2000}.
    We shall show that it is unstable if $\gamma > 1$ and Lyapunov stable if $\gamma < -3$ as well; see \Cref{prop:eqilateral_triangle_with_center}.
  }
  \label{fig:equilateral_with_center}
\end{figure}

In order to demonstrate applications of \Cref{thm:Energy-Casimir}, let us first consider the relative equilibria with four vortices at the vertices and the center of the unit equilateral triangle (i.e., each side is the unit length) as shown in \Cref{fig:equilateral_with_center}.
At the vertices of the equilateral triangle are three vortices of circulation 1, whereas at the center is one with an arbitrary circulation $\gamma \in \R \backslash \{-3\}$ so that the total circulation $\Gamma = \sum_{i=1}^{4} \Gamma_{i}$ is non-zero.

\citet{CaSc2000} claims to have proved the Lyapunov stability for $-1/2 < \gamma < 1$ using a different method that is more specialized for regular polygon configuration of $N$ vortices with a central vortex.
However, this result is contested by \citet{KuOsSo2016} (see also \cite{KuOs2021,KuOsSo2016,KuOsSo2025}); in particular it seems stable for $\gamma = -0.6$ in their numerical experiment.
It is also proved in \cite{KuOsSo2016} that it is stable if $\gamma < -3$ or $0 < \gamma < 1$.
We shall illustrate \Cref{thm:Energy-Casimir} by reproducing this result.

Since this is the case with $N = 4$ and $\Gamma \neq 0$, the Lie algebra $\vK$ is the set of $3 \times 3$ skew-Hermitian matrices, and so the matrix $\mu$ takes the form
\begin{equation*}
  \mu = \rmi
  \begin{bmatrix}
    \mu_{1} & \mu_{12} & \mu_{13} \\
    \mu_{12}^{*} & \mu_{2} & \mu_{23} \\
    \mu_{13}^{*} & \mu_{23}^{*} & \mu_{3}
  \end{bmatrix}
  = \rmi
  \begin{bmatrix}
    \mu_{1} & \mu_{4} + \rmi \mu_{5} & \mu_{6} + \rmi\mu_{7} \\
    \mu_{4} - \rmi \mu_{5} & \mu_{2} & \mu_{8} + \rmi\mu_{9} \\
    \mu_{6} - \rmi\mu_{7} & \mu_{8} - \rmi\mu_{9} & \mu_{3}
  \end{bmatrix},
\end{equation*}
and so we shall use the coordinates $(\mu_{1}, \dots, \mu_{9}) \in \R^{9}$.

\begin{proposition}[Stability of equilateral triangle with center~\cite{KuOsSo2016}]
  \label[proposition]{prop:eqilateral_triangle_with_center}
  The fixed point $\mu_{0} \in \vK$ of the Lie--Poisson relative dynamics~\eqref{eq:Lie-Poisson} corresponding to the relative equilibrium of an equilateral triangle with center (see \Cref{fig:equilateral_with_center}) is Lyapunov stable to $F^{-1}(0)$-preserving perturbations if $\gamma < -3$ or $0 < \gamma < 1$ and linearly unstable if $\gamma > 1$.
\end{proposition}
\begin{proof}[Proof of \Cref{prop:eqilateral_triangle_with_center}]
  The fixed point is given by
  \begin{equation*}
    \mu_{0} = \parentheses{ 1, 1, 1, -\frac{1}{2}, -\frac{\sqrt{3}}{2}, -\frac{1}{2}, \frac{\sqrt{3}}{2}, -\frac{1}{2}, -\frac{\sqrt{3}}{2} },
  \end{equation*}
  and the linearization of the Lie--Poisson equation~\eqref{eq:Lie-Poisson} at $\mu_{0}$ gives a linear system in $\R^{9}$ with $9\times9$ matrix with eigenvalues
  \begin{equation*}
    \pm\frac{1}{\pi}\rmi,\
    \pm\frac{1}{2\pi}\sqrt{\gamma - 1},\
    0,
  \end{equation*}
  where the algebraic multiplicity of each of $\pm\frac{1}{2\pi}\sqrt{\gamma - 1}$ is 2, and that of $0$ is $3$.
  Hence we see that $\mu_{0}$ is linearly unstable if $\gamma > 1$.

  Let us now apply \Cref{thm:Energy-Casimir} to this case.
  The constraints are given by the zero level set of
  \begin{equation*}
    F\colon \nonzerovK \to \R^{4};
    \qquad
    F(\mu) =
    \begin{bmatrix}
      F_{1}(\mu) \smallskip\\
      F_{2}(\mu) \smallskip\\
      \Re F_{12}(\mu) \smallskip\\
      \Im F_{12}(\mu)
    \end{bmatrix}
    =
    \begin{bmatrix}
      \mu_{1} \mu_{2} - \mu_{4}^{2} - \mu_{5}^{2} \smallskip\\
      \mu_{2} \mu_{3} - \mu_{8}^{2} - \mu_{9}^{2} \smallskip\\
      -\mu_{2} \mu_{6} + \mu_{4} \mu_{8} - \mu_{5} \mu_{9} \smallskip\\
      -\mu_{2} \mu_{7} + \mu_{5} \mu_{8} + \mu_{4} \mu_{9}
    \end{bmatrix}.
  \end{equation*}
  One can then show that, among the Casimirs $\{ C_{j} \}_{j=1}^{3}$ shown in \eqref{eq:Casimirs}, only one (and any one) of them is independent of $F$ at $\mu_{0}$.
  Hence we shall take the simplest one
  \begin{equation*}
    C_{1}(\mu) = \frac{1}{\gamma + 3}\parentheses{ (\gamma+2) \mu_{1} + (\gamma+2) \mu_{2} + 2\mu_{3} + \gamma \mu_{3} - 2\mu_{4} - 2\mu_{6} - 2\mu_{8} },
  \end{equation*}
  and define
  \begin{equation}
    \label{eq:f-equilateral_with_center}
    f(\mu) \defeq 4\pi a_{0} h(\mu) + a_{1} C_{1}(\mu) + b_{1} F_{1}(\mu) + b_{2} F_{2}(\mu) + c_{12} \Re F_{12}(\mu) + d_{12} \Im F_{12}(\mu).
  \end{equation}
  It is a straightforward computation to show that $Df(\mu_{0}) = 0$ if
  \begin{equation*}
    (a_{1}, b_{1}, b_{2}, c_{12}, d_{12}) = a_{0} \parentheses{
      \gamma + 1,\,
      \frac{2\gamma}{3(\gamma + 3)},\,
      \frac{2\gamma}{3(\gamma + 3)},\,
      -\frac{4\gamma}{3(\gamma + 3)},\,
      0
    }.
  \end{equation*}
  Now define the submanifold
  \begin{equation*}
    M \defeq F^{-1}(0) \cap C_{1}^{-1}(C_{1}(\mu_{0})) \subset \vK \cong \R^{9}.
  \end{equation*}
  Then the following vectors $\{ v_{i} \}_{i=1}^{4}$ give a basis for the tangent space $T_{\mu_{0}}M$:
  \begin{gather*}
    v_{1} = \sqrt{3}e_{1} - \sqrt{3}e_{3} - e_{5},
    \quad
    v_{2} = e_{1} - e_{3} - e_{4} + e_{8},
    \\
    v_{3} = -\sqrt{3}e_{2} + \sqrt{3}e_{3} + e_{5} + e_{7},
    \quad
    v_{4} = e_{2} - e_{3} - e_{4} + e_{6}
  \end{gather*}
  using the standard basis $\{ e_{i} \}_{i=1}^{9}$ for $\R^{9}$.
  Defining a $4 \times 4$ matrix $\mathcal{H}$ by setting $\mathcal{H}_{ij} \defeq v_{i}^{T} D^{2}f(\mu_{0}) v_{j}$, its leading principal minors are
  \begin{gather*}
    d_{1} = a_{0} \frac{2(9\gamma^{2} + 20\gamma + 3)}{3(\gamma + 3)},
    \qquad
    d_{2} = -a_{0}^{2} \frac{16\gamma (\gamma - 1)}{3 (\gamma +3)},
    \\
    d_{3} = -a_{0}^{3} \frac{8 \gamma (\gamma-1) (9\gamma^{2} + 20\gamma + 3)}{3 (\gamma+3)^{2}},
    \qquad
    d_{4} = a_{0}^{4} \frac{16 (\gamma -1)^{2} \gamma ^{2}}{(\gamma +3)^{2}}.
  \end{gather*}
  The polynomial $9\gamma^{2} + 20\gamma + 3$ in $d_{1}$ and $d_{3}$ have roots $(-10 \pm \sqrt{73})/9 \simeq -2.06, -0.162$.
  Thus one sees the following two cases under which $d_{i} > 0$ for every $i \in \{1, 2, 3, 4\}$:
  (i)~$a_{0} > 0$ and $0 < \gamma < 1$; (ii)~$a_{0} < 0$ and $\gamma < -3$ (see \Cref{fig:dets-EquilateralWithCenter}).
  \begin{figure}[htbp]
    \centering
    \begin{subcaptionblock}[c]{0.475\textwidth}
      \centering
      \includegraphics[width=\textwidth]{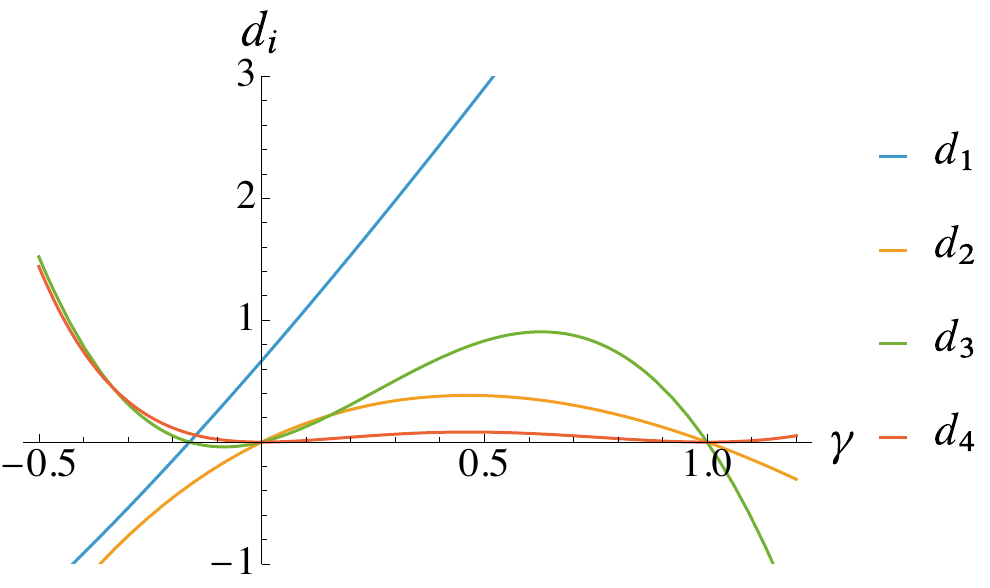}
      \caption{$a_{0} = 1$}
      \label{fig:dets-EquilateralWithCenter1}
    \end{subcaptionblock}
    \hfill
    \begin{subcaptionblock}[c]{0.475\textwidth}
      \centering
      \includegraphics[width=\textwidth]{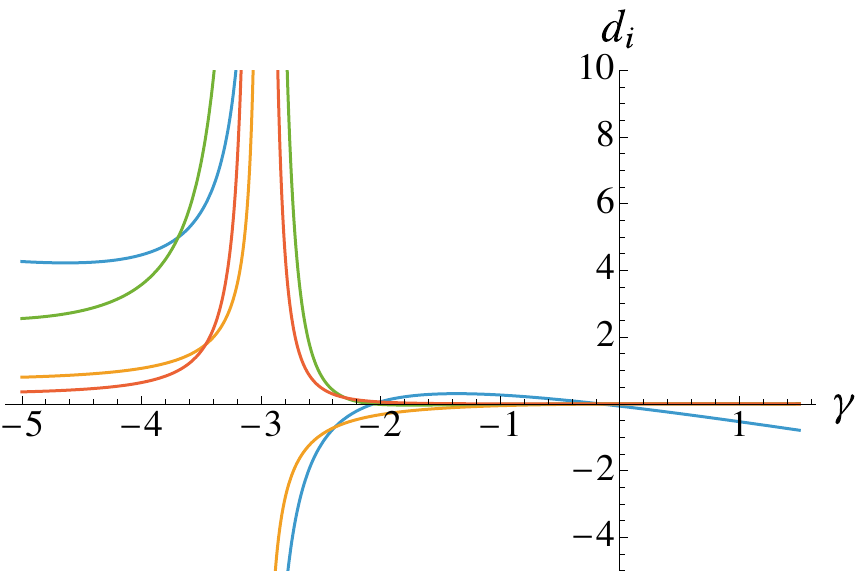}
      \caption{$a_{0} = -1/10$}
      \label{fig:dets-EquilateralWithCenter2}
    \end{subcaptionblock}
    \caption{Leading principal minors $d_{i}$ with $i \in \{1, 2, 3, 4\}$ of $\mathcal{H}$ with \subref{fig:dets-EquilateralWithCenter1}~$a_{0} = 1$ and \subref{fig:dets-EquilateralWithCenter2}~$a_{0} = -1/10$.}
    \label{fig:dets-EquilateralWithCenter}
  \end{figure}
  
  Therefore, by \Cref{thm:Energy-Casimir}, the fixed point $\mu_{0}$ is Lyapunov stable if $\gamma < -3$ or $0 < \gamma < 1$.
\end{proof}

\subsection{Stability of Rhombus Equilibria}
\begin{figure}[htbp]
  \centering
  \includegraphics[width=.65\textwidth]{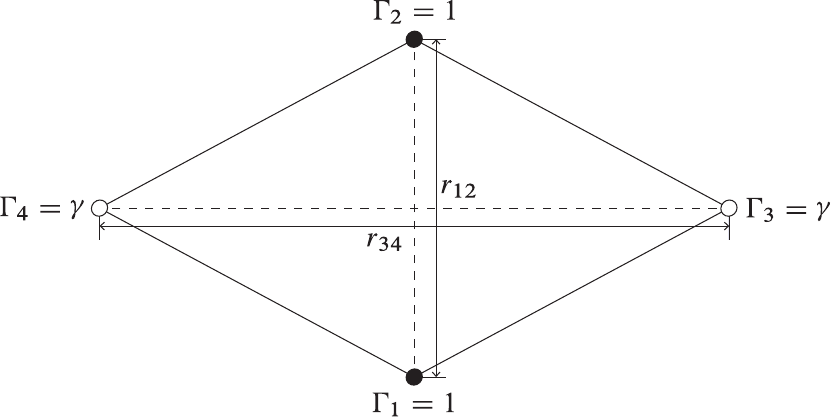}
  \caption{
    The rhombus shown above give a relative equilibria (called the Rhombus-A equilibria in \cite{Ro2013,HaRoSa2014}) if the inter-vortex distances $r_{12}$ and $r_{34}$ satisfy the relationship given in \eqref{eq:alpha-gamma}.
    It is known to be linearly stable if $-2 + \sqrt{3} < \gamma < 0$, Lyapunov stable for $0 < \gamma \le 1$, and unstable if $-1 < \gamma < -2 + \sqrt{3}$~\cite{Ro2013}.
    We shall show that it is Lyapunov stable for $-2 + \sqrt{3} < \gamma < 0$, extending the domain of $\gamma$ for Lyapunov stability.
  }
  \label{fig:Rhombus-A}
\end{figure}

Let us also consider a one-parameter family of rhombus relative equilibria shown in \Cref{fig:Rhombus-A}.
This is one of the families of relative equilibria found in \citet[Theorem~7.7]{HaRoSa2014} consisting of two pairs of vortices of same strengths: vortices 1 and 2 with circulation 1 and vortices 3 and 4 with $\gamma$, where $\gamma \in (-1,1]$; note that we assume $\gamma \neq -1$ to avoid the vanishing total circulation, i.e., $\Gamma = \sum_{i=1}^{4} \Gamma_{i} \neq 0$.
The cases with $\gamma < -1$ or $\gamma > 1$ can be rescaled to yield another Rhombus-A equilibrium or one of the other family of rhombus equilibria called Rhombus-B in \cite{HaRoSa2014}, which is known to be always unstable~\cite{Ro2013}.

\begin{figure}
  \centering
  \includegraphics[width=.55\linewidth]{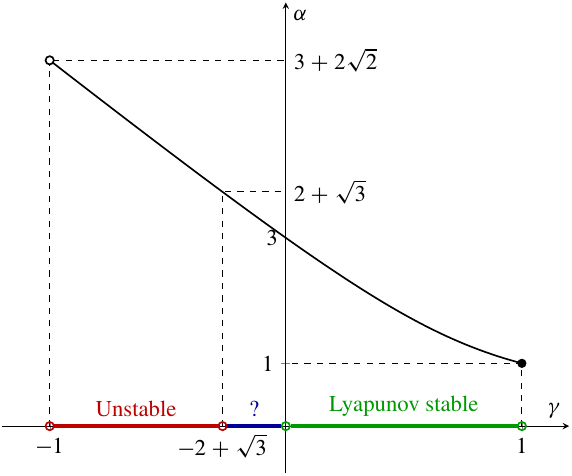}
  \caption{Relationship between parameter $\alpha$ and circulation $\gamma$; see \eqref{eq:alpha-gamma}. $-1 < \gamma \le 1$ corresponds to $1 \le \alpha < 3 + 2\sqrt{2}$. Rhombus-A relative equilibria have been known to be unstable for $-1 < \gamma < -2 + \sqrt{3}$ and Lyapunov stable for $0 < \gamma \le 1$, but only linearly stable if $-2 + \sqrt{3} < \gamma < 0$~\cite{Ro2013}. We shall prove in \Cref{prop:Rhombus-A} that they are Lyapunov stable for $-2 + \sqrt{3} < \gamma < 0$ as well.}
  \label{fig:alpha-gamma}
\end{figure}

It was shown in \cite[Theorem~7.7]{HaRoSa2014} that it gives a family of relative equilibria if the distances $r_{ij} \defeq |q_{i} - q_{j}|$ between the vortices satisfy
\begin{equation}
  \label{eq:alpha-gamma}
  \parentheses{ \frac{r_{34}}{r_{12}} }^{2} = \frac{1}{2} \parentheses{ 3(1 - \gamma) + \sqrt{ 9(1 - \gamma)^{2} + 4\gamma } } \eqdef \alpha.
\end{equation}
We shall use $\alpha$ instead of $\gamma$ as the parameter of the problem as computations are much simpler that way; see \Cref{fig:alpha-gamma} for the plot of $\alpha$ against $\gamma$ defined above for $-1 < \gamma \le 1$.

\citet{Ro2013} proved that the Rhombus-A relative equilibria are linearly stable if $-2 + \sqrt{3} < \gamma < 0$ and Lyapunov stable for $0 < \gamma \le 1$, whereas it is unstable if $-1 < \gamma < -2 + \sqrt{3}$.
We note that one can consider the rhombus equilibria as a special case of the double alternate rings considered in \citet{Ha1931}, and the stability condition there gives the same linear stability condition.

As is well known, a linearly stable equilibrium of a Hamiltonian system can be unstable in the sense of Lyapunov; see, e.g., \citet{Ch1928} and also \citet[Chapters~13]{MeHaOf2008} for its exposition.
Therefore, the only problem of interest for us is the Lyapunov stability for $-2 + \sqrt{3} < \gamma < 0$, which corresponds to $3 < \alpha < 2 + \sqrt{3}$; see \Cref{fig:alpha-gamma}.

\begin{proposition}[Stability of Rhombus-A]
  \label[proposition]{prop:Rhombus-A}
  The fixed point $\mu_{0} \in \vK$ of the Lie--Poisson relative dynamics~\eqref{eq:Lie-Poisson} corresponding to the relative equilibrium of Rhombus-A (see \Cref{fig:Rhombus-A}) is Lyapunov stable to $F^{-1}(0)$-preserving perturbations if $-2 + \sqrt{3} < \gamma < 0$ or $0 < \gamma < 1$, and linearly unstable if $-1 < \gamma < -2 + \sqrt{3}$.
\end{proposition}
\begin{proof}
  Using the parameter $\alpha$ defined in \eqref{eq:alpha-gamma}, the fixed point is given by
  \begin{equation*}
    \mu_{0} = \parentheses{ \alpha + 1,\, \alpha + 1,\, 4\alpha,\, \alpha - 1,\, -2\sqrt{\alpha},\, 2\alpha,\, -2\sqrt{\alpha},\, 2\alpha,\, 2\sqrt{\alpha}\, },
  \end{equation*}
  and the linearization of the Lie--Poisson equation~\eqref{eq:Lie-Poisson} at $\mu_{0}$ gives a linear system in $\R^{9}$ with $9\times9$ matrix with eigenvalues
  \begin{equation*}
    \pm\frac{\alpha - 3}{4\pi(\alpha + 1)}\rmi,\
    \pm\frac{ \sqrt{F_{1}(\alpha)} }{\pi(\alpha + 1)(3\alpha - 1)},\
    \pm\frac{ \sqrt{F_{2}(\alpha)} }{4\pi(\alpha + 1)^{2}(3\alpha - 1)},\
    0,
  \end{equation*}
  where
  \begin{equation*}
    F_{1}(\alpha) \defeq \parentheses{ \alpha^{2} - 4\alpha + 1 } \parentheses{ 3\alpha^{2} - 2\alpha + 3 }
    \begin{cases}
      < 0 & 1 \le \alpha < 2 + \sqrt{3}, \\
      = 0 & \alpha = 2 + \sqrt{3}, \\
      > 0 & 2 + \sqrt{3} < \alpha < 3 + 2\sqrt{2},
    \end{cases}
  \end{equation*}
  and
  \begin{equation*}
    F_{2}(\alpha) \defeq \parentheses{ 3\alpha^{3} - 15\alpha^{2} + 41\alpha - 5 } \parentheses{ 5\alpha^{3}- 41\alpha^{2} + 15\alpha - 3 } < 0
    \quad
    \text{for}
    \quad
    1 \le \alpha < 3 + 2\sqrt{2}.
  \end{equation*}
  Therefore, we see that $\mu_{0}$ is linearly unstable if $2 + \sqrt{3} < \alpha < 3 + 2\sqrt{2}$, reproducing \cite{Ro2013}.

  The Lyapunov stability analysis works the same way as in the proof of \Cref{prop:eqilateral_triangle_with_center}.
  Using the same form of Lyapunov function $f$ from \eqref{eq:f-equilateral_with_center}, we have $Df(\mu_{0}) = 0$ if
  \begin{equation*}
    a_{1} = -\frac{\alpha^{2} - 14\alpha + 1}{6\alpha^{2} + 4 \alpha - 2}\, a_{0},
    \quad
    b_{1} = b_{2} = \frac{c_{12}}{2} = -\frac{\alpha(\alpha - 3)^{2}}{4(\alpha + 1)^{2}(\alpha^{2} - 6\alpha + 1)}\, a_{0},
    \quad
    d_{12} = 0.
  \end{equation*}
  One may take a basis for the tangent space $T_{\mu_{0}}M$ as
  \begin{gather*}
    v_{1} = 2\sqrt{\alpha}\, e_{1} + \sqrt{\alpha}\, e_{4} - e_{5} + 2\sqrt{\alpha}\, e_{6}, \qquad
    v_{2} = 2\sqrt{\alpha}\, e_{2} + \sqrt{\alpha}\, e_{4} - e_{5} + 2\sqrt{\alpha}\, e_{8}, \\
    v_{3} = -e_{1} + e_{2} + \sqrt{\alpha}\,e_{7} + \sqrt{\alpha}\,e_{9}, \\
    v_{4} = -(\alpha + 1)(3\alpha - 1) e_{1} + (2\alpha^{3} - 9\alpha^{2} + 4\alpha - 1) e_{2} + 4\alpha(3\alpha - 1) e_{3} \\
    \qquad - \alpha(\alpha^{2} - 1) e_{4} - \sqrt{\alpha}(\alpha + 1)(\alpha^{2} - 4\alpha + 1) e_{5} + 2\sqrt{\alpha} (\alpha - 1)^{3} e_{9}.
  \end{gather*}
  Let us define a $4 \times 4$ matrix $\mathcal{H}$ by setting 
  \begin{equation*}
    \mathcal{H}_{ij} \defeq \frac{2(\alpha + 1)^{2} (3\alpha - 1) (-\alpha^{2} + 6\alpha - 1)}{\alpha(\alpha - 3)}\, v_{i}^{T} D^{2}f(\mu_{0}) v_{j},
  \end{equation*}
  where we placed the rational function in $\alpha$ (which is positive for $\alpha > 3$) as a multiplicative factor in order to simplify the calculations to follow.
  Then the leading principal minors $\{ d_{i} \}_{i=1}^{4}$ of $\mathcal{H}$ satisfy
  \begin{gather*}
    d_{1} = a_{0}(13\alpha^{3} - 89\alpha^{2} + 23\alpha - 3),
    \qquad
    d_{2} = 32 a_{0}^{2} (-\alpha^{2} + 6\alpha - 1) (5\alpha^{3} - 41\alpha^{2} + 15\alpha - 3), \\
    \frac{d_{3}}{d_{2}} = -a_{0}(3\alpha^{3} - 15\alpha^{2} + 41\alpha - 5), \\
    \frac{d_{4}}{d_{3}} = -2a_{0}\, \alpha (\alpha + 1)^{2} (-\alpha^{2} + 6\alpha - 1) (-\alpha^{2} + 4\alpha - 1) (3\alpha^{2} - 2\alpha + 3).
  \end{gather*}
  Pick an arbitrary $a_{0} < 0$.
  Then, one can easily verify that $d_{i} > 0$ for $1 \le i \le 4$ when $1 < \alpha < 3$ as well as $3 < \alpha < 2 + \sqrt{3}$, which correspond to $0 < \gamma < 1$ (again reproducing \cite{Ro2013}) and $-2 + \sqrt{3} < \gamma < 0$, respectively; see \eqref{fig:alpha-gamma}.
\end{proof}

\appendix
\numberwithin{equation}{section}

\section{Proof of \Cref{prop:constraint_functions}~(\ref{prop:constraint_functions-submersion})}
\label{sec:submersion}
We would like to prove that $F$ defined in \eqref{eq:F} is a submersion, i.e., the functions $\{ F_{i} \}_{i=1}^{N-2} \cup \{ F_{ij} \}_{1\le i < j \le N-2}$ are independent on $\nonzerovK$ in the sense that their differentials evaluated at every point in $\nonzerovK$ are linearly independent.

Let us first rewrite the non-diagonal entries of $\mu$ as
\begin{equation*}
  \mu_{ij} = x_{ij} + \rmi y_{ij}
  \quad\text{with}\quad
  x_{ij}, y_{ij} \in \R
  \quad\text{for}\quad
  1 \le i < j \le N-2
\end{equation*}
so that we may identify $\vK$ with $\R^{(N-1)^{2}}$ as follows:
\begin{align*}
  \vK \cong \R^{(N-1)^{2}} &\cong \R^{N-1} \times \R^{(N-1)(N-2)} \\
                           &= \bigl\{ \bigl( \mu_{1}, \dots, \mu_{N}, (x_{12}, y_{12}), \dots, (x_{N-2,N-1}, y_{N-2,N-1}) \bigr) \bigr\}.
\end{align*}
For brevity we set $n \defeq N - 1$, and shall prove the independence by induction on $n$.
For $n = 2$, we have only $F_{1}(\mu) = \mu_{1} \mu_{2} - (x_{12}^{2} + y_{12}^{2})$, and the differential $\d{F_{1}}(\mu)$ is clearly non-zero when $\mu \in \nonzerovK$.
Suppose that the claimed independence holds for $n - 1$.
We may write the matrix $-\rmi\mu$ as follows:
\begin{equation*}
  -\rmi\mu = 
  \begin{bNiceMatrix}
    \mu_{1}        & \mu_{12} & \Cdots          & \mu_{1,n-1} & \mu_{1,n}\\
    \mu_{12}^{*}    & \Ddots  & \Ddots          & \vdots  & \vdots \\
    \vdots         & \Ddots  & \Ddots          & \mu_{n-2,n-1} & \mu_{n-2,n} \\
    \mu_{1,n-1}^{*} & \Cdots  & \mu_{n-2,n-1}^{*} & \mu_{n-1} & \mu_{n-1,n} \\
    \mu_{1,n}^{*} & \Cdots  & \mu_{n-2,n}^{*} & \mu_{n-1,n}^{*} & \mu_{n}
  \end{bNiceMatrix}
\end{equation*}
We shall show that the functions $\{ F_{i,n-1} \}_{i=1}^{n-2} \cup \{ F_{n-1} \}$ defined by the $2 \times 2$ determinants sweeping the last two columns of the above matrix are independent of those preceding ones $\{ F_{i,j} \}_{1\le i < j \le n-2} \cup \{ F_{i} \}_{i=1}^{n-2}$ defined for the $(n-1) \times (n-1)$ submatrix obtained by removing the last row and the last column.

We see that, for $1 \le i \le n-4$, using simplified subscripts for brevity,
\begin{align*}
  F_{i,n-1}(\mu) =
  \begin{vmatrix}
    \mu_{i,n-1} & \mu_{i,n} \\
    \mu_{i+1,n-1} & \mu_{i+1,n}
  \end{vmatrix}
  &=
    \begin{vmatrix}
      x_{1} + \rmi y_{1} & x_{4} + \rmi y_{4} \\
      x_{2} + \rmi y_{2} & x_{5} + \rmi y_{5}
    \end{vmatrix},
  \\
  F_{i+1,n-1}(\mu) =
  \begin{vmatrix}
    \mu_{i+1,n-1} & \mu_{i+1,n} \\
    \mu_{i+2,n-1} & \mu_{i+2,n}
  \end{vmatrix}
  &=
    \begin{vmatrix}
      x_{2} + \rmi y_{2} & x_{5} + \rmi y_{5} \\
      x_{3} + \rmi y_{3} & x_{6} + \rmi y_{6}
    \end{vmatrix},
\end{align*}
and we see that, written as row vectors,
\begin{equation*}
  \begin{bmatrix}
    \d(\Re F_{i,n-1})(\mu) \\
    \d(\Im F_{i,n-1})(\mu) \\
    \d(\Re F_{i+1,n-1})(\mu) \\
    \d(\Im F_{i+1,n-1})(\mu) 
  \end{bmatrix}
  =
  \begin{bNiceMatrix}
    \Block{4-4}<\Large>{*}
    & & & & -x_{2} &  y_{2} & x_{1} & -y_{1} & 0 & 0 \\
    & & & & -y_{2} & -x_{2} & y_{1} &  x_{1} & 0 & 0 \\
    & & & &  0 & 0 & -x_{3} & y_{3} & x_{2} & -y_{2} \\
    & & & &  0 & 0 & -y_{3} & -x_{3} & y_{2} & x_{2}
  \end{bNiceMatrix},
\end{equation*}
where we have explicitly written only the last six columns, which represent coefficients of
\begin{equation*}
  \{ \d{x}_{4}, \d{y}_{4}, \d{x}_{5}, \d{y}_{5}, \d{x}_{6}, \d{y}_{6} \}.
\end{equation*}
Since $\{ (x_{i}, y_{i}) \}_{i=4}^{6}$ are the variables in $\mu_{i,n}$, $\mu_{i+1,n}$, and $\mu_{i+2,n}$ from the last column of $\mu$, the preceding determinants $\{ F_{i,j} \}_{1\le i < j \le n-2} \cup \{ F_{i} \}_{i=1}^{n-2}$ do not depend on them, hence their differentials have vanishing coefficients for $\{ \d{x}_{i}, \d{y}_{i} \}_{i=4}^{6}$.
Now, since $\mu \in \nonzerovK$, $x_{i} + \rmi y_{i} \neq 0$ for $i = 1, 2, 3$, and so the $4 \times 6$ submatrix explicitly shown above is full rank, implying that the four differentials $\d(\Re F_{i,n-1})(\mu)$, $\d(\Im F_{i,n-1})(\mu)$, $\d(\Re F_{i+1,n-1})(\mu)$, and $\d(\Im F_{i+1,n-1})(\mu)$ are independent of each other as well as independent of the differentials of the preceding ones.

One can apply a similar argument about the pair $F_{n-3,n-1}$ and $F_{n-2,n-1}$, where $F_{n-2,n-1}$ looks slight different:
\begin{equation*}
  F_{n-2,n-1}(\mu) =
  \begin{vmatrix}
    \mu_{n-2,n-1} & \mu_{n-2,n} \\
    \mu_{n-1} & \mu_{n-1,n}
  \end{vmatrix}
  =
  \begin{vmatrix}
    x_{2} + \rmi y_{2} & x_{4} + \rmi y_{4} \\
    x_{3} & x_{5} + \rmi y_{5}
  \end{vmatrix}.
\end{equation*}
Hence this is a special case of the above argument with $y_{3} = 0$ and $x_{3} \neq 0$ by assumption, and so again implies the same independence property.

It now remains to show that the differential of the last determinant $F_{n-1}$ is independent of the rest.
In fact,
\begin{equation*}
  F_{n-1}(\mu) =
  \begin{vmatrix}
    \mu_{n-1} & \mu_{n-1,n} \\
    \mu_{n-1,n}^{*} & \mu_{n}
  \end{vmatrix}
  = \mu_{n-1} \mu_{n} - |\mu_{n-1,n}|^{2},
\end{equation*}
and thus the term proportional to $\d{\mu}_{n}$ in $\d{F_{n-1}}(\mu)$ is $\mu_{n-1}\d{\mu}_{n}$ where $\mu_{n-1} \neq 0$.
But then all the other determinants are independent of $\mu_{n}$.
Hence $\d{F_{n-1}}(\mu)$ is independent of all the other differentials.

\section*{Acknowledgments}
This work was supported by NSF grant DMS-2006736.
I would like to thank Holger Dullin, Luis Garc\'ia-Naranjo, Richard Montgomery, and the reviewers for their helpful comments and suggesting some of the references.

\bibliographystyle{plainnat}
\bibliography{Point_Vortex_Stability-Plane}

\begin{thebibliography}{42}
\providecommand{\natexlab}[1]{#1}
\providecommand{\url}[1]{\texttt{#1}}
\expandafter\ifx\csname urlstyle\endcsname\relax
  \providecommand{\doi}[1]{doi: #1}\else
  \providecommand{\doi}{doi: \begingroup \urlstyle{rm}\Url}\fi

\bibitem[Aeyels(1992)]{Ae1992}
D.~Aeyels.
\newblock On stabilization by means of the {E}nergy--{C}asimir method.
\newblock \emph{Systems \& Control Letters}, 18\penalty0 (5):\penalty0
  325--328, 1992.

\bibitem[Albouy and Dullin(2020)]{AlDu2020}
A.~Albouy and H.~R. Dullin.
\newblock Relative equilibria of the 3-body problem in $\mathbb{R}^4$.
\newblock \emph{Journal of Geometric Mechanics}, 12\penalty0 (3):\penalty0
  323--341 EP --, 2020.

\bibitem[Aref(2007)]{Ar2007}
H.~Aref.
\newblock Point vortex dynamics: A classical mathematics playground.
\newblock \emph{Journal of Mathematical Physics}, 48\penalty0 (6):\penalty0
  065401, 2007.

\bibitem[Aref(1979)]{Ar1979}
H.~Aref.
\newblock Motion of three vortices.
\newblock \emph{The Physics of Fluids}, 22\penalty0 (3):\penalty0 393--400,
  1979.

\bibitem[Arunasalam et~al.(2015)Arunasalam, Dullin, and Nguyen]{ArDuNg2015}
S.~Arunasalam, H.~R. Dullin, and D.~M.~H. Nguyen.
\newblock {The Lie--Poisson structure of the symmetry reduced regularized
  $n$-body problem}.
\newblock \emph{Journal of Physics A: Mathematical and Theoretical},
  48\penalty0 (6):\penalty0 065202, 2015.

\bibitem[Barry and Hoyer-Leitzel(2016)]{BaHo2016}
A.~M. Barry and A.~Hoyer-Leitzel.
\newblock Existence, stability, and symmetry of relative equilibria with a
  dominant vortex.
\newblock \emph{SIAM Journal on Applied Dynamical Systems}, 15\penalty0
  (4):\penalty0 1783--1805, 2016.

\bibitem[Bolsinov et~al.(1999)Bolsinov, Borisov, and Mamaev]{BoBoMa1999}
A.~V. Bolsinov, A.~V. Borisov, and I.~S. Mamaev.
\newblock Lie algebras in vortex dynamics and celestial mechanics---{IV}.
\newblock \emph{Regular and Chaotic Dynamics}, 4\penalty0 (1):\penalty0 23--50,
  1999.

\bibitem[Borisov and Pavlov(1998)]{BoPa1998}
A.~V. Borisov and A.~E. Pavlov.
\newblock Dynamics and statics of vortices on a plane and a sphere---{I}.
\newblock \emph{Regular and Chaotic Dynamics}, 3\penalty0 (1):\penalty0 28--38,
  1998.

\bibitem[Borisov and Mamaev(2015)]{BoMa2015}
A.~V. Borisov and I.~S. Mamaev.
\newblock Symmetries and reduction in nonholonomic mechanics.
\newblock \emph{Regular and Chaotic Dynamics}, 20\penalty0 (5):\penalty0
  553--604, 2015.

\bibitem[Cabral and Schmidt(2000)]{CaSc2000}
H.~E. Cabral and D.~S. Schmidt.
\newblock Stability of relative equilibria in the problem of {$N+1$} vortices.
\newblock \emph{SIAM Journal on Mathematical Analysis}, 31\penalty0
  (2):\penalty0 231--250, 2000.

\bibitem[Cherry(1928)]{Ch1928}
T.~M. Cherry.
\newblock On periodic solutions of hamiltonian systems of differential
  equations.
\newblock \emph{Philosophical Transactions of the Royal Society of London.
  Series A, Containing Papers of a Mathematical or Physical Character},
  227:\penalty0 137--221, 1928.

\bibitem[Chorin and Marsden(1993)]{ChMa1993}
A.~J. Chorin and J.~E. Marsden.
\newblock \emph{A Mathematical Introduction to Fluid Mechanics}, volume~4 of
  \emph{Texts in Applied Mathematics}.
\newblock Springer, 1993.

\bibitem[Dullin and Montgomery()]{DuMo2025}
H.~Dullin and R.~Montgomery.
\newblock The {$N$}-body problem on coadjoint orbits.
\newblock \emph{https://arxiv.org/abs/2504.01272}.

\bibitem[Dullin(2013)]{Du2013}
H.~R. Dullin.
\newblock {The Lie--Poisson structure of the reduced $n$-body problem}.
\newblock \emph{Nonlinearity}, 26\penalty0 (6):\penalty0 1565, 2013.

\bibitem[Dullin and Scheurle(2020)]{DuSc2020}
H.~R. Dullin and J.~Scheurle.
\newblock Symmetry reduction of the 3-body problem in $\mathbb{R}^{4}$.
\newblock \emph{Journal of Geometric Mechanics}, 12\penalty0 (3):\penalty0
  377--394, 2020.

\bibitem[Hampton et~al.(2014)Hampton, Roberts, and Santoprete]{HaRoSa2014}
M.~Hampton, G.~E. Roberts, and M.~Santoprete.
\newblock Relative equilibria in the four-vortex problem with two pairs of
  equal vorticities.
\newblock \emph{Journal of Nonlinear Science}, 24\penalty0 (1):\penalty0
  39--92, 2014.

\bibitem[Havelock(1931)]{Ha1931}
T.~Havelock.
\newblock {LII. The stability of motion of rectilinear vortices in ring
  formation}.
\newblock \emph{The London, Edinburgh, and Dublin Philosophical Magazine and
  Journal of Science}, 11\penalty0 (70):\penalty0 617--633, 1931.

\bibitem[Kurakin and Ostrovskaya(2021)]{KuOs2021}
L.~G. Kurakin and I.~V. Ostrovskaya.
\newblock Resonances in the stability problem of a point vortex quadrupole on a
  plane.
\newblock \emph{Regular and Chaotic Dynamics}, 26\penalty0 (5):\penalty0
  526--542, 2021.

\bibitem[Kurakin et~al.(2016)Kurakin, Ostrovskaya, and Sokolovskiy]{KuOsSo2016}
L.~G. Kurakin, I.~V. Ostrovskaya, and M.~A. Sokolovskiy.
\newblock {On the stability of discrete tripole, quadrupole, Thomson'vortex
  triangle and square in a two-layer/homogeneous rotating fluid}.
\newblock \emph{Regular and Chaotic Dynamics}, 21\penalty0 (3):\penalty0
  291--334, 2016.

\bibitem[Kurakin et~al.(2025)Kurakin, Ostrovskaya, and Sokolovskiy]{KuOsSo2025}
L.~G. Kurakin, I.~V. Ostrovskaya, and M.~A. Sokolovskiy.
\newblock On the stability of discrete $n+1$ vortices in a two-layer rotating
  fluid: The cases $n=4,5,6$.
\newblock \emph{Regular and Chaotic Dynamics}, 30\penalty0 (3):\penalty0
  325---353, 2025.

\bibitem[Lerman and Singer(1998)]{LeSi1998}
E.~Lerman and S.~F. Singer.
\newblock Stability and persistence of relative equilibria at singular values
  of the moment map.
\newblock \emph{Nonlinearity}, 11\penalty0 (6):\penalty0 1637, 1998.

\bibitem[Marsden and Ratiu(1999)]{MaRa1999}
J.~E. Marsden and T.~S. Ratiu.
\newblock \emph{Introduction to Mechanics and Symmetry}.
\newblock Springer, 1999.

\bibitem[Marsden et~al.(2007)Marsden, Misiolek, Ortega, Perlmutter, and
  Ratiu]{MaMiOrPeRa2007}
J.~E. Marsden, G.~Misiolek, J.~P. Ortega, M.~Perlmutter, and T.~S. Ratiu.
\newblock \emph{Hamiltonian Reduction by Stages}.
\newblock Springer, 2007.

\bibitem[McLachlan et~al.(2016)McLachlan, Modin, and Verdier]{McMoVe2016b}
R.~I. McLachlan, K.~Modin, and O.~Verdier.
\newblock Symmetry reduction for central force problems.
\newblock \emph{European Journal of Physics}, 37\penalty0 (5):\penalty0 055003,
  2016.

\bibitem[Menezes and Roberts(2018)]{MeRo2018}
B.~Menezes and G.~E. Roberts.
\newblock Existence and stability of four-vortex collinear relative equilibria
  with three equal vorticities.
\newblock \emph{SIAM Journal on Applied Dynamical Systems}, 17\penalty0
  (1):\penalty0 1023--1051, 2018.

\bibitem[Meyer et~al.(2008)Meyer, Hall, and Offin]{MeHaOf2008}
K.~Meyer, G.~Hall, and D.~Offin.
\newblock \emph{{Introduction to Hamiltonian Dynamical Systems and the $N$-Body
  Problem}}.
\newblock Springer, 2nd edition, 2008.

\bibitem[Montaldi(2000)]{Mo2000}
J.~Montaldi.
\newblock {Relative equilibria and conserved quantities in symmetric
  Hamiltonian systems}.
\newblock In R.~Kaiser and J.~Montaldi, editors, \emph{Peyresq Lectures On
  Nonlinear Phenomena}. World Scientific Publishing Company, 2000.

\bibitem[Montaldi(1997)]{Mo1997}
J.~Montaldi.
\newblock Persistence and stability of relative equilibria.
\newblock \emph{Nonlinearity}, 10\penalty0 (2):\penalty0 449, 1997.

\bibitem[Montaldi and Rodr{\'\i}guez-Olmos(2011)]{MoRo2011}
J.~Montaldi and M.~Rodr{\'\i}guez-Olmos.
\newblock {On the stability of Hamiltonian relative equilibria with non-trivial
  isotropy}.
\newblock \emph{Nonlinearity}, 24\penalty0 (10):\penalty0 2777, 2011.

\bibitem[Montgomery(2015)]{Mo2015}
R.~Montgomery.
\newblock The three-body problem and the shape sphere.
\newblock \emph{American Mathematical Monthly}, 122\penalty0 (4):\penalty0
  299--321, 2015.

\bibitem[Newton(2001)]{Ne2001}
P.~K. Newton.
\newblock \emph{The $N$-vortex problem}.
\newblock Springer, New York, 2001.

\bibitem[Ohsawa(2019)]{Oh2019d}
T.~Ohsawa.
\newblock Symplectic reduction and the {L}ie--{P}oisson shape dynamics of {$N$}
  point vortices on the plane.
\newblock \emph{Nonlinearity}, 32\penalty0 (10):\penalty0 3820--3842, 2019.

\bibitem[Ohsawa(2023)]{Oh2023a}
T.~Ohsawa.
\newblock Shape dynamics of {$N$} point vortices on the sphere.
\newblock \emph{Nonlinearity}, 36\penalty0 (2):\penalty0 1000--1028, 2023.

\bibitem[Ortega and Ratiu(1999)]{OrRa1999}
J.-P. Ortega and T.~S. Ratiu.
\newblock {Stability of Hamiltonian relative equilibria}.
\newblock \emph{Nonlinearity}, 12\penalty0 (3):\penalty0 693, 1999.

\bibitem[Patrick(1992)]{Pa1992}
G.~W. Patrick.
\newblock {Relative equilibria in Hamiltonian systems: The dynamic
  interpretation of nonlinear stability on a reduced phase space}.
\newblock \emph{Journal of Geometry and Physics}, 9\penalty0 (2):\penalty0
  111--119, 1992.
\newblock ISSN 0393-0440.

\bibitem[Patrick et~al.(2004)Patrick, Roberts, and Wulff]{PaRoWu2004}
G.~W. Patrick, M.~Roberts, and C.~Wulff.
\newblock Stability of poisson equilibria and hamiltonian relative equilibria
  by energy methods.
\newblock \emph{Archive for Rational Mechanics and Analysis}, 174\penalty0
  (3):\penalty0 301--344, 2004.

\bibitem[P\'erez-Chavela et~al.(2015)P\'erez-Chavela, Santoprete, and
  Tamayo]{PeSaTa2015}
E.~P\'erez-Chavela, M.~Santoprete, and C.~Tamayo.
\newblock Symmetric relative equilibria in the four-vortex problem with three
  equal vorticities.
\newblock \emph{Dynamics of Continuous, Discrete and Impulsive Systems Series
  A: Mathematical Analysis}, 22:\penalty0 189--209, 2015.

\bibitem[Roberts(2013)]{Ro2013}
G.~E. Roberts.
\newblock Stability of relative equilibria in the planar $n$-vortex problem.
\newblock \emph{SIAM Journal on Applied Dynamical Systems}, 12\penalty0
  (2):\penalty0 1114--1134, 2013.

\bibitem[Roberts(2018)]{Ro2018}
G.~E. Roberts.
\newblock Morse theory and relative equilibria in the planar $n$-vortex
  problem.
\newblock \emph{Archive for Rational Mechanics and Analysis}, 228\penalty0
  (1):\penalty0 209--236, 2018.

\bibitem[Souriau(1997)]{So1997}
J.~Souriau.
\newblock \emph{Structure of Dynamical Systems: A Symplectic View of Physics}.
\newblock Progress in Mathematics. Birkh{\"a}user Boston, 1997.

\bibitem[Sundaram(1996)]{Su1996}
R.~Sundaram.
\newblock \emph{A First Course in Optimization Theory}.
\newblock Cambridge University Press, 1996.

\bibitem[Synge(1949)]{Sy1949}
J.~L. Synge.
\newblock On the motion of three vortices.
\newblock \emph{Canadian Journal of Mathematics}, 1\penalty0 (3):\penalty0
  257--270, 1949.

\end{thebibliography}

\end{document}